\documentclass[11pt,twoside]{article}
\usepackage{amsfonts}
\usepackage{amsmath}
\usepackage{amsthm}
\usepackage{graphicx}

\usepackage{enumerate}

\usepackage[urlcolor=blue]{hyperref}

\newcommand{\bi}{\begin{itemize}}\newcommand{\ei}{\end{itemize}}
\newcommand{\be}{\begin{equation}}\newcommand{\ee}{\end{equation}}
\newcommand{\bee}{\begin{enumerate}}\newcommand{\eee}{\end{enumerate}}
\newcommand{\bea}{\begin{eqnarray}}\newcommand{\eea}{\end{eqnarray}}
\newcommand{\beas}{\begin{eqnarray*}}\newcommand{\eeas}{\end{eqnarray*}}
\newcommand{\bc}{\begin{center}}\newcommand{\ec}{\end{center}}

\def\qed{
\vbox{\hrule\hbox{\vrule\hbox to 5pt{\vbox to 8pt{\vfil}\hfil}\vrule}\hrule}}
\def\qedsmall{
\vbox{\hrule\hbox{\vrule\hbox to 4pt{\vbox to 4pt{\vfil}\hfil}\vrule}\hrule}}

\def\endproof{\unskip \nobreak \hskip0pt plus 1fill \qquad \qed}

\def\9{\"u9}
\newtheoremstyle{custom}{3pt}{3pt}{}{}{\bfseries}{:}{.5em}{}
\theoremstyle{custom}

\numberwithin{equation}{section}
\numberwithin{figure}{section}
\numberwithin{table}{section}

\def\R{\mathbb{R}}

\def\N{\mathbb{N}}
\def\Z{\mathbb{Z}}
\def\X{\mathbb{X}}

\def\U{\mathbb{U}}

\def\KK{{\cal K}}

\def\UU{{\cal U}}

\def\QQ{{\cal Q}}

\def\eps{\varepsilon}

\def\beq{\begin{equation}}
\def\eeq{\end{equation}}
\def\bea{\begin{eqnarray}}
\def\eea{\end{eqnarray}}
\def\beaa{\begin{eqnarray*}}
\def\eeaa{\end{eqnarray*}}

\usepackage{color}
\def\LG{\textcolor{black}}

\title{Dissipativity and optimal control}
\author{Lars Gr\"une\\Chair of Applied Mathematics, University of 
Bayreuth, Germany\\
	lars.gruene@uni-bayreuth.de}

\begin{document}
\maketitle

\pagestyle{myheadings}
\thispagestyle{plain}
\markboth{LARS GR\"UNE}{DISSIPATIVITY AND OPTIMAL CONTROL}

\section{Introduction}

The close link between dissipativity and optimal control is already apparent in Jan C. Willems' first papers on the subject. Particularly, the paper \cite{Will71}, which appeared one year before his famous dissipativity papers \cite{Will72a,Will72b}, already contains a lot of insight on this connection for linear quadratic optimal control problems. \LG{Around the same time, the role of a passivity-like variant of dissipativity for inverse optimal control of nonlinear systems was described by Moylan and Anderson \cite{MoyA73}, a result that triggered a significant amount of research with a peak of activities in the 1990s and monographs such as those by Freeman and Kokotovi\'{c} \cite{FreK96} or by Sepulchre, Jankovic and Kokotovi\'{c} \cite{SeJK97}.}
In recent years, research on the link \LG{between dissipativity and optimal control} has been revived with a particular focus on \LG{the turnpike phenomenon} and applications in model predictive control (MPC). 

This development was initiated about ten years ago in the paper \cite{DiAR10} by Diehl, Amrit, and Rawlings, in which 
\LG{strong duality} was used in order to construct a Lyapunov function for the closed-loop solution resulting from an economic MPC scheme. Soon after it was realized in \cite{AnAR12} that \LG{strong duality is nothing else than strict dissipativity with a linear storage function, and that} linearity of the storage function is not really needed for this result, meaning that the same Lyapunov function construction can be carried out for general nonlinear strict dissipativity. In both cases, the link between dissipativity and optimal control lies in the fact that the running cost (or the stage cost in discrete time) serves as the supply rate in the dissipativity formulation. 

While these results use strict dissipativity of the optimal control problem within the MPC scheme in a formal way in order to obtain stability results via an appropriate Lyapunov function, the effect of strict dissipativity on optimal trajectories can also be analyzed in a more geometrical fashion. More precisely, generalizing techniques that were already around in the 1990s (see, e.g., \cite[Theorem 4.2]{CaHL91}), it was observed in \cite[Theorems 5.3 and 5.6]{Grue13} that strict dissipativity plus a suitable controllability property is sufficient for the occurrence of the so-called turnpike property. This property, first observed in mathematical economy in the works of Ramsey and von Neumann in the 1920s and 1930s \cite{Rams28,vNeu38}, and first called this way by Dorfman, Samuelson and Solow in the 1950s \cite{DoSS58}, describes the fact that optimal trajectories most of the time stay close to an optimal equilibrium. 

Based on this observation, in \cite{GruS14,FauB15,GruP15a,Grue16} (see also \cite[Chapter 7]{GruP17} and \cite{FaGM18}) the stability results from \cite{AnAR12} could be extended to larger classes of MPC schemes and, in addition, non-averaged or transient approximate optimality of the MPC closed-loop could be established. Motivated by these new applications, the relation of strict dissipativity to classical notions of detectability for linear and nonlinear systems was also clarified \cite{GruG18,GruG21,HoeG19}. This paper surveys these recent developments and some of Willems' early results.

\section{Optimal control problems and (strict) dissipativity}

\subsection{Optimal control problems}

We consider optimal control problems either in continuous time
\be \mbox{minimize } J_T(x_0,u) = \int_0^T \ell(x(t),u(t)) dt \label{eq:ctoc} \ee
with respect to $u\in\UU$, $u(t)\in\U$ and $x(t)\in \X$ for all $t\in[0,T]$, $T\in\R_{>0}$, where $\UU$ is an appropriate space of functions, $\X\subset\R^n$ and $\U\subset\R^m$ are the sets of \emph{admissible states} and \emph{admissible control inputs}, respectively, and
\be \dot x(t) = f(x(t),u(t)), \qquad x(0)= x_0, \label{eq:ctsys}\ee
or in discrete time 
\be \mbox{minimize } J_T(x_0,u) = \sum_{\LG{t=0}}^{T-1} \ell(x(t),u(t)) \label{eq:dtoc} \ee
with respect to $u\in\UU$, $u(t)\in\U$ and $x(t)\in \X$ for all $t=0,\ldots,T$, $T\in\N$, where $\UU$ is an appropriate space of sequences, $\X$ and $\U$ are as above, and
\be x(t+1) = f(x(t),u(t)), \qquad x(0)= x_0, \label{eq:dtsys}\ee
which we briefly write as $x^+=f(x,u)$. Here $f:\R^n\times \R^m \to \R^n$ is either the vector field in continuous time or the iteration map in discrete time, while $\ell:\R^n\times \R^m \, \LG{\to \R}$ is called the running cost in continuous time and the stage cost in discrete time. In order to unify the notation, we use the symbol $[a,b]$ both in continuous and discrete time. In continuous time it denotes the usual closed interval $\{t\in\R\,|\, a\le t\le b\}$, while in discrete time it denotes $\{t\in\Z\,|\, a\le t\le b\}$, where $\Z$ is the set of integers. For simplicity of exposition we limit \LG{ourselves} to finite dimensional state space but we mention that some of the results discussed in this paper are also available in infinite dimensional settings. 
Given an initial value $x_0\in\X$, we denote the set of control functions $u\in\UU$ for which $x(t)\in\X$ \LG{ and $u(t)\in\U$} holds for all $t\in[0,T]$ by $\UU(x_0,T)$. The \emph{optimal value function} is then defined as
\[ V_T(x_0) := \inf_{u\in\UU(x_0,T)} J_T(x_0,u) \]
and a control $u^*\in\UU(x_0,T)$ with corresponding trajectory $x^*(\cdot)$ is called \emph{optimal control} for initial condition $x_0$ and time horizon $T$ if 
\[ J_T(x_0,u^*) = V_T(x_0).\] 
The corresponding trajectory $x^*(\cdot)$ is then called an \emph{optimal trajectory}.

\subsection{Dissipativity and optimal control}

Dissipativity in the sense of Willems as we use it in this paper involves an abstract notion of energy that is stored in the system. For each admissible state $x\in\R^n$, we denote the energy in the system by $\lambda(x)$. The function $\lambda:\R^n\to\R$ is called \emph{storage function} and it is usually assumed that $\lambda$ is bounded from below, in order to avoid that an infinite amount of energy can be extracted from the system. In continuous time, dissipativity then demands that there exists another function, the so called \emph{supply rate} $s:\R^n\times \R^m \to \R$, such that the inequality
\be \lambda(x(\tau)) \le \lambda(x_0) + \int_0^\tau s(x(t),u(t))dt \label{eq:ctdiss}\ee
holds for all $\tau\ge 0$, all control functions $u\in\UU(x_0,\tau)$ and all initial conditions $x_0\in\X$. As in the optimal control problem \eqref{eq:ctoc}, $x(t)$ is supposed to satisfy \eqref{eq:ctsys}. In discrete time, the demanded inequality completely analogously reads
\be \lambda(x(\tau)) \le \lambda(x_0) + \sum_{t=0}^{\tau-1} s(x(t),u(t)),  \label{eq:dtdiss}\ee
with $x(t)$ satisfying \eqref{eq:dtsys}. It appears that this discrete-time variant of \eqref{eq:ctdiss} was first used by Byrnes and Lin in \cite{ByrL94}. 

The interpretation of these inequalities is as follows: They demand that the energy $\lambda(x(\tau))$ in the system after a certain time $\tau$ is not larger than the initial energy $\lambda(x_0)$ plus the integral or sum over the supplied energy, expressed at each time instant by $s(x(t),\LG{u(t)})$. \LG{In this interpretation, if} $s(x(t),\LG{u(t)})$ \LG{is negative, then} negative energy is supplied, \LG{which means that} energy is extracted from the system.

In continuous time, if $\lambda$ is continuously differentiable, \LG{using that because of the chain rule we have
\[ \frac{d}{dt} \lambda(x(t)) = D\lambda(x(t))\dot x(t) = D\lambda(x(t))f(x(t),u(t)),\]
}
inequality \eqref{eq:ctdiss} can equivalently be rewritten in infinitesimal form
\be D\lambda(x)f(x,u) \le s(x,u), \label{eq:ctdinf}\ee
while in discrete time equivalently to \eqref{eq:dtdiss} one can use the one-step form
\be \lambda(x^+) \le \lambda(x) + s(x,u),  \label{eq:dtdone}\ee
where we use the common brief notation $x^+ = f(x,u)$. \LG{In order to establish these equivalences, the inequalities \eqref{eq:ctdinf} and \eqref{eq:dtdone} need to hold for all pairs $(x,u)=(x(t),u(t))$ that appear in \eqref{eq:ctdiss} or \eqref{eq:dtdiss}, respectively.}

This dissipativity concept can be related to the optimal control problem from the previous section by setting $s(x,u) := \ell(x,u)$ for the running or stage cost $\ell$ from \eqref{eq:ctoc} and \eqref{eq:dtoc}, respectively. It appears that this connection was first made in Willems' paper \cite{Will71}, which remarkably appeared in the year \emph{before} his seminal papers \cite{Will72a,Will72b}, which introduced and discussed dissipativity in a comprehensive way. More precisely, in \cite{Will71} continuous-time, linear, controllable dynamics $f(x,u) = Ax + Bu$ and quadratic running cost $\ell(x,u) = x^TQx + 2u^TCx + u^TRu$ without any definiteness assumption on $Q$ and $R$ are studied. The paper provides necessary and sufficient conditions for the existence of a storage function $\lambda$ in terms of the finiteness of optimal value functions of certain related optimal control problems. These characterizations led to the concepts of \emph{available storage} and \emph{required supply} that are described in the sidebar on \nameref{sidebar:avstorage}. 
What is important here \LG{is} that these characterizations involve constraints on the asymptotic behavior of the optimal control problems under consideration. This already indicates the fundamental connection between dissipativity and the long-term behavior of optimal trajectories\LG{, an inspiring research topic until today, as recent publications such as \cite{FauK20} show.}

\subsection{\LG{Dissipativity and inverse optimal control}}

\LG{Passivity is a particular case of dissipativity with supply rate $s(x,u)=u^Ty = u^Th(x)$ or variants thereof, where $y=h(x)$ is the output of the system. As stated in the introduction of his paper \cite{Will72a}, Willems introduced dissipativity as a ``generalization of the concept of passivity''. The link between passivity and inverse optimal control goes back to \cite{MoyA73}. Here we describe it following the presentation in \cite[Section 3.3]{SeJK97}.}

\LG{
To this end, consider the continuous-time functional \eqref{eq:ctoc} on the infinite time horizon, i.e., with $T=\infty$ and with a running cost $\ell$ of the form $\ell(x,u) = \hat \ell(x) + \|u\|^2$. Consider, moreover, dynamics \eqref{eq:ctsys} with $f$ of the form $f(x,u) = f_1(x) + f_2(x)u$. Then there exists a cost function $\hat\ell\ge 0$ such that the feedback control $u=-k(x)$ is optimal and stabilizing if and only if the closed-loop system
\begin{eqnarray*} \dot x(t) & = & f_1(x(t)) + f_2(x(t))u(t) \\
y(t) & = & k(x)
\end{eqnarray*}
is zero-state detectable and dissipative with supply rate $s(x,u) = u^Tk(x) + \frac12 k(x)^Tk(x)$ and continuously differentiable storage function $\lambda$. Here, zero-state detectability means that every solution with $u\equiv 0$ and $y\equiv 0$ behaves like a solution of an asymptotically stable system.
}

\LG{The term ``inverse'' in inverse optimality stems from the fact that the cost function $\hat \ell$ is not known but, instead, only implicitly defined by $f$ and $k$ (for details we refer to the references just cited). Rather than specifying $\hat \ell$ and deriving $k$, here one specifies or designs $k$ and derives $\hat \ell$. One may then ask why inverse optimality is an important property, given that one has only limited control about what is actually optimized by the feedback law $k$. The answer is that optimal feedback laws have a range of beneficial properties, above all related to robustness, and we refer to \cite{FreK96} and \cite{SeJK97} and the references therein for details.
}

\LG{The fact that $k$ is a stabilizing feedback law is crucial for this result. Hence, again, there is a 
connection between dissipativity and the long-term behavior of optimal trajectories. Summarizing, the two classical lines of research 
\begin{itemize}
\item revealed that optimal control problems can indicate whether dissipativity holds, e.g., via the available storage or the required supply, \textit{or}
\item showed that dissipativity conditions on stabilizing feedback laws imply that they solve an optimal control problem for a suitable cost function.
\end{itemize}
The recent renaissance of the investigation of the link between dissipativity and optimal control, which we will describe in the remainder of this paper, looks at a third aspect. Namely, its goal is to
\begin{itemize}
\item use dissipativity in order to make statements about the long-time behavior of optimal and near-optimal trajectories.
\end{itemize}
}

Before we turn to the description of this recent development, we introduce a stricter variant of the dissipativity property.

\subsection{Strict dissipativity}

Strict dissipativity is nothing but dissipativity with storage function 
\[ s(x,u) - \alpha(\|x\LG{-x^e}\|) \]
in place of $s(x,u)$, where $\alpha\in\KK_\infty$ with
\[ \KK_\infty:=\{\alpha:[0,\infty)\to[0,\infty)\,|\, \alpha \mbox{ continuous, strictly \LG{increasing}, unbounded, and }\alpha(0)=0\}\]
and $x^e$ is an equilibrium of the control system, i.e.\ there is a control value $u^e$ with $f(x^e,u^e)=0$ in continuous time and $f(x^e,u^e)=x^e$ in discrete time. Written explicitly, the corresponding dissipation inequalities read 
\be \lambda(x(\LG{\tau})) \le \lambda(x_0) + \int_0^{\LG{\tau}} s(x(\LG{t}),u(\LG{t})) - \alpha(\|x(\LG{t})-x^e\|) d\LG{t} \label{eq:ctsdiss}\ee
and 
\be \lambda(x(\LG{\tau})) \le \lambda(x_0) + \sum_{\LG{t=0}}^{\LG{\tau-1}} s(x(\LG{t}),u(\LG{t})) - \alpha(\|x(\LG{t})-x^e\|)\LG{.}  \label{eq:dtsdiss}\ee
The \LG{infinitesimal or one-step} versions \eqref{eq:ctdinf} and \eqref{eq:dtdone} of these inequalities then become
\be D\lambda(x)f(x,u) \le s(x,u)-\alpha(\|x-x^e\|), \label{eq:ctdinfs}\ee
and
\be \lambda(x^+) \le \lambda(x) + s(x,u)-\alpha(\|x-x^e\|),  \label{eq:dtdones}\ee
respectively. 

For some applications, it is necessary to use $\alpha(\|(x,u)-(x^e,u^e)\|)$ in place of $\alpha(\|x-x^e\|)$. In this case, one speaks about \emph{strict $(x,u)$-dissipativity}. 

One easily sees that \eqref{eq:ctsdiss} and \eqref{eq:dtsdiss} are more demanding than \eqref{eq:ctdiss} and \eqref{eq:dtdiss}, respectively, since they do not only demand that the energy difference $\lambda(x(t))-\lambda(x_0)$ is bounded by the integral or sum over the supplied energy $s$, but that actually some of the energy is dissipated if the system is not in the equilibrium $x^e$, and the amount of dissipated energy increases the further away the state $x(\tau)$ is from $x^e$. 

At a first glance it seems that the difference between strict and non-strict dissipativity is only quantitative. Indeed, it follows readily from the definition that if the system is dissipative with supply $s(x,u)$ then it is strictly dissipative with supply $s(x,u)+\alpha(\|x\|)$. However, if the supply function $s$ is not merely a parameter we can play with but a function that results from some modelling procedure, then there is not only a quantitative but also a qualitative difference between strict and non-strict dissipativity. This in particular applies for the setting in which $s$ is derived from the running or stage cost $\ell$ of an optimal control problem, i.e., when $s=\ell$. 

When using strict dissipativity for the analysis of the control system behavior, it is often necessary to demand $s(x^e,u^e)=0$, as we will see in the next section. This is in general a restrictive condition, but it is actually not restrictive in case that $s$ is derived from the optimal control problem. This is because the optimal trajectories for the cost $\ell(x,u)$ are the same as for the cost $\ell(x,u)-\ell(x^e,u^e)$. It is thus convenient and not restrictive to use the supply rate
\be s(x,u) = \ell(x,u) - \ell(x^e,u^e). \label{eq:optsupply}\ee

\section{The turnpike property}

The turnpike property demands that optimal trajectories (and in some variants also near-optimal trajectories) stay in the vicinity of an equilibrium $x^e$ most of the time. Of course, this statement needs to be made mathematically precise in order to be able to analyze it mathematically. The meaning of ``most of the time'' is that the amount of time the trajectory spends outside a given neighborhood of the equilibrium $x^e$ is bounded, and the bound is independent of the optimization horizon and of the initial condition, at least for initial conditions which itself are contained in a bounded set. \LG{In this paper we focus on the turnpike property for finite-horizon optimal control problems and we refer to \cite{GrKW17} for an infinite-horizon version and the relation between the finite-horizon and infinite-horizon case and to \cite{FauK20} for the impact of dissipativity on infinite-horizon optimal control.}

\begin{figure}[htb]
\includegraphics[width=14cm]{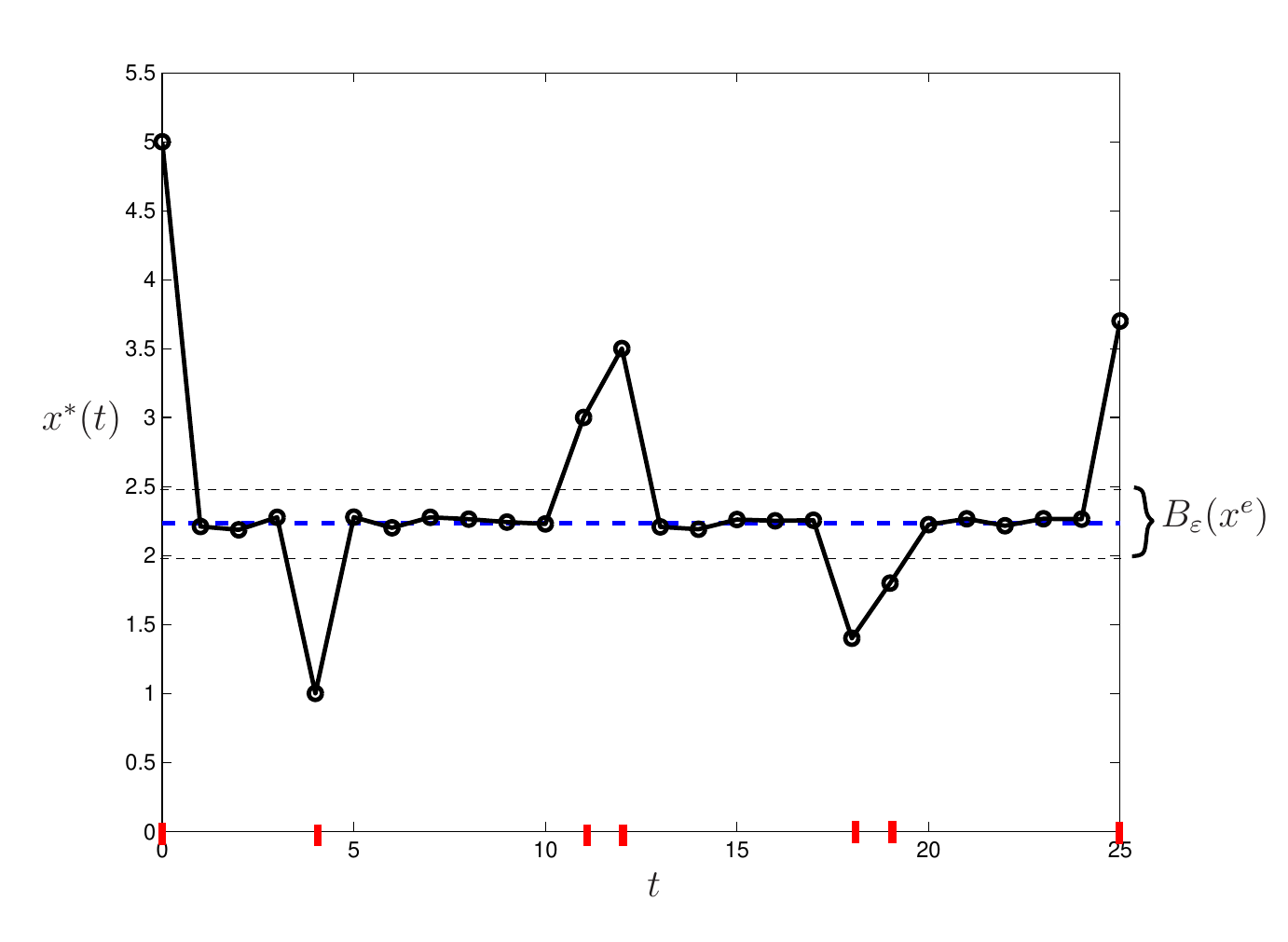}
\caption{Illustration of the turnpike property in discrete time. The set $\QQ$ contains $7$ time instants, which are marked with short red vertical lines on the $t$-axis. The equilibrium $x^e$ is represented by the dashed blue line.\label{fig:turnpike_illustration}}
\end{figure}

Formally, we say that an optimal control problem has the \emph{turnpike property}, if there is an equilibrium $x^e$, such that for any $\eps>0$ and $K>0$ there is a $C>0$ such that for all time horizons $T>0$ and all optimal trajectories $x^*(\cdot)$ with $x^*(0)\in B_K(x^e)\cap\X$ the set of times
\[ \QQ := \{ t\in [0,T]\,|\, x^*(t)\not\in B_\eps(x^e)\} \]
satisfies $|\QQ|\le C$. Here $|\QQ|$ is the Lebesgue measure of $\QQ$ in continuous time or the number of elements contained in $\QQ$ in discrete time, and $B_K(x^e)$ and $B_\eps(x^e)$ denote the balls around $x^e$ with radii $K$ and $\eps$, respectively. 
Figure \ref{fig:turnpike_illustration} shows a sketch of a trajectory exhibiting the turnpike phenomenon. Here the set $\QQ$ contains $7$ time instants, which are marked with short red vertical lines on the $t$-axis. The equilibrium $x^e$ is represented by the dashed blue line.

\begin{figure}[htb]
\includegraphics[width=7cm]{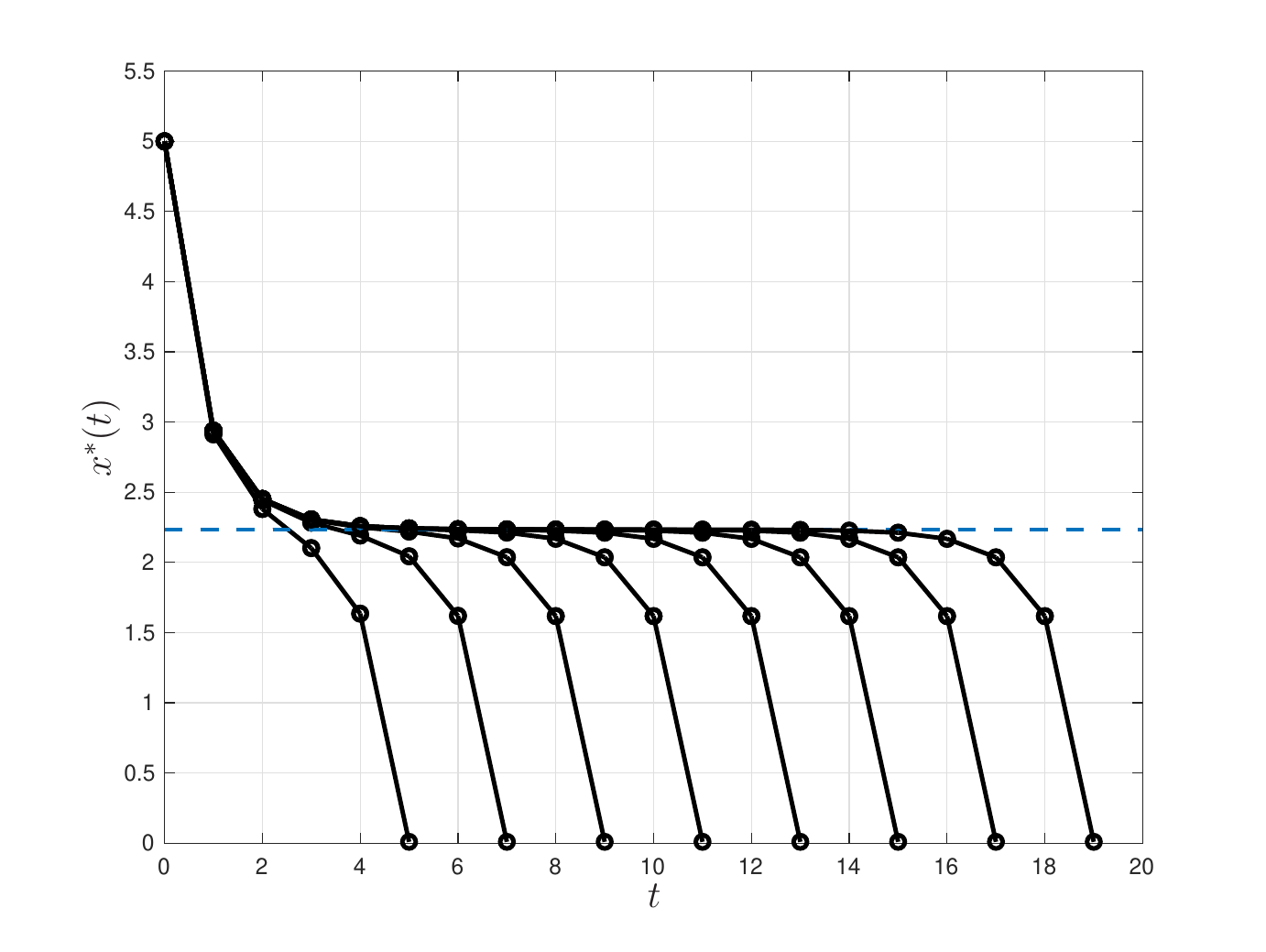}
\includegraphics[width=7cm]{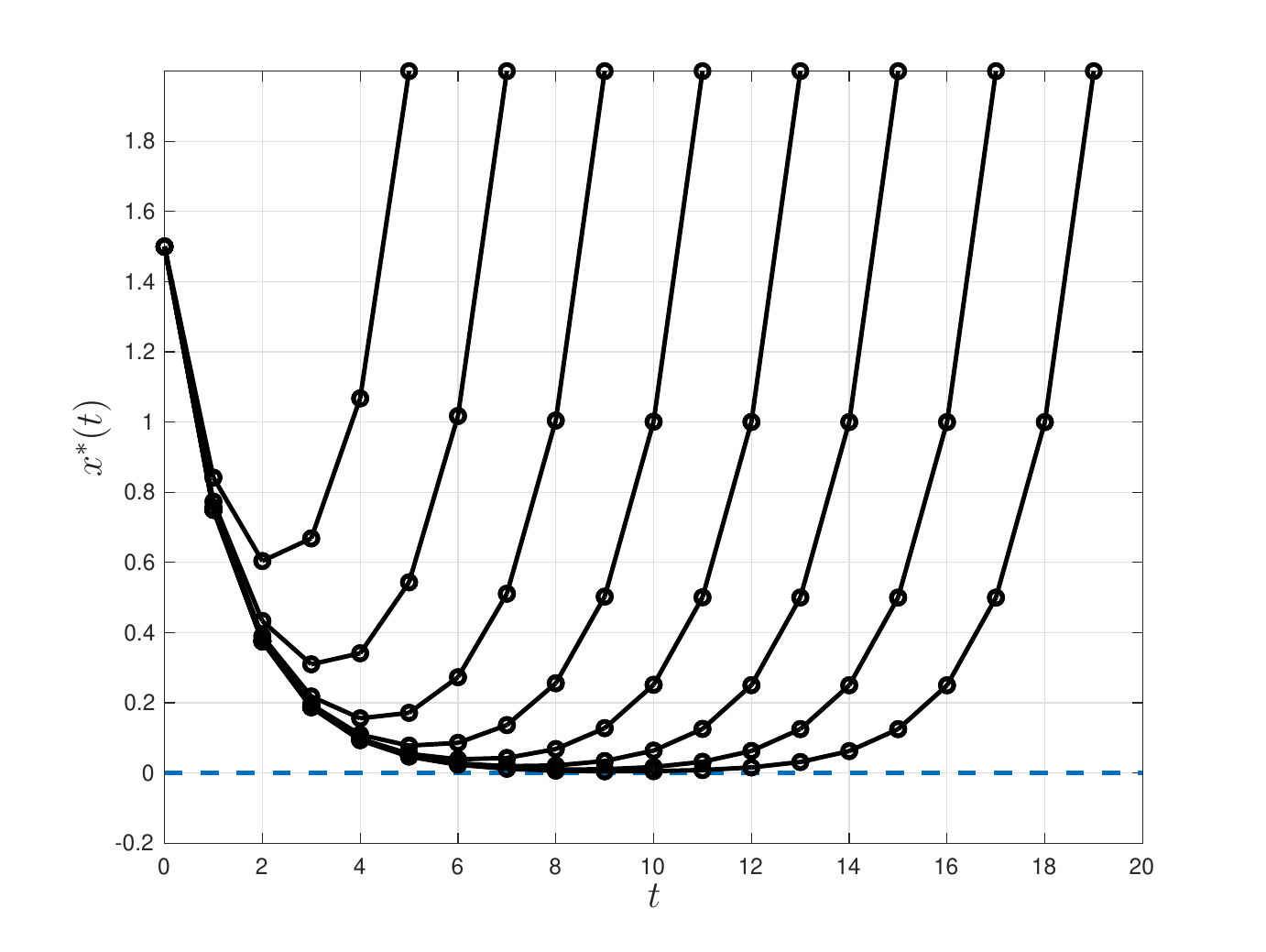}
\caption{Optimal trajectories for the examples \eqref{eq:brock} and \eqref{eq:invariance} for time horizon length $T=5,7,\ldots,19$.\label{fig:turnpike_examples}}
\end{figure}

Figure \ref{fig:turnpike_examples} shows optimal trajectories for varying time horizons $T=5,7,\ldots,19$ for two simple one-dimensional discrete-time optimal control problems. The first problem is given by 
\be x^+ = u, \quad \ell(x,u) = -\log(5x^{0.34} - u), \quad \X=[0,10], \quad \U = [0.01,10] \label{eq:brock}\ee
and the second by 
\be x^+ = 2x + u, \quad \ell(x,u) = u^2, \quad \X=[-2,2], \quad \U = [-3,3].\label{eq:invariance}\ee
The model \eqref{eq:brock} is an optimal investment problem from \cite{BroM72}, in which $x$ denotes the investment in a company and $5x^{0.34}$ is the return from this investment (including the investment itself) after one time period. As $x^+=u$ is the investment in the next time period, $5x^{0.34} - u$ is the amount of money that can be used for consumption in the current time period and the optimal control problem thus models the maximization of the sum of the logarithmic utility function $\log(5x^{0.34} - u)$ over the time periods. Clearly, it is optimal to spend all the available money until the end of the time horizon $T$, which is why all optimal trajectories end at $x=0$. However, in between the trajectories spend most of the time near the equilibrium $x^e \approx 2.2344$, i.e.\ it exhibits the turnpike property. We will explain later why this is the case.

The task modelled in the second example is to keep the state of the system in $\X$ with as little quadratic control effort $\ell(x,u)=u^2$ as possible. In the long run it is beneficial to stay near $x^e=0$, because the control effort for staying in $\X$ is very small near this equilibrium. Thus, it makes sense for the optimal trajectories to stay near $x^e=0$ for most of the time. At the end of the horizon, however, it is beneficial to turn off the control completely, because this entirely reduces the control effort to $0$ and does not violate the state constraint provided the control is turned off sufficiently late.

In both examples in Figure \ref{fig:turnpike_examples} it is clearly visible that the number of states outside a neighborhood of the respective equilibrium (indicated by the dashed blue line) remains constant with increasing horizon length, which is exactly what the turnpike property demands.

These examples show that the turnpike property occurs already in very simple problems. It is, however, known that it also occurs in many much more complicated problems, including optimal control problems governed by partial differential equations; \cite{PorZ13,GuTZ16,Zuaz17,Porr18,TrZZ18,GrSS19,BreP20,GrSS20} is just a small selection of recent references where this is investigated. 

Interestingly, the two examples exhibit a particular form of the turnpike property, in which the exceptional points, i.e., the points corresponding to the times in $\QQ$, lie only at the beginning and at the end of the time interval. One speaks of the \emph{approaching arc} at the beginning and the \emph{leaving arc} at the end of the time horizon. While the general definition of the turnpike property allows for excursions from $x^e$ also in the middle of the time interval (as indicated in the sketch in Figure \ref{fig:turnpike_illustration}), these do not appear in the two examples. The reason for this will be explained in the next section. 

The optimal trajectories in Figure \ref{fig:turnpike_examples} nicely illustrate the source of the name ``turnpike property'', which was coined in \cite{DoSS58}: the behavior of the optimal trajectory is similar to a car driving from the initial to the end point, where the equilibrium---i.e., the dashed blue line---plays the role of a highway, or turnpike, as highways are called in parts of the USA. If the time is sufficiently large, it pays off to first go to the highway (even if this is associated with some additional cost), stay there for most of the time and then leave the highway at the end. There may be different reasons for the occurence of the leaving arc. We may impose a constraint $x(N)=\hat x$ for the terminal state that forces the trajectory to leave the turnpike, or it may simply be beneficial to leave the turnpike because this reduces the cost of the overall trajectory. As we did not impose any terminal constraints, the latter must be the case in our examples, and it is actually easy to convince yourself that this is the case.

Besides the definition given above, there are a number of alternative definitions of the turnpike property. A widely used variant is the \emph{exponential turnpike property}, which demands an inequality of the type 
\[ \|x^*(t) - x^e\| \le C(e^{-\sigma t} + e^{-\sigma(T-t)}) \]
for all $t\in[0,T]$ with constants $C,\sigma>0$. On easily sees that this property implies the turnpike property above, since for each $\eps>0$ there is $\tau>0$ with $C(e^{-\sigma_1 t} + e^{-\sigma_2(T-t)}) < \eps$ for all $t\in[\tau,T-\tau]$ holds, regardless of how large $T$ is. 
Another variant demands the turnpike behavior not only for the optimal trajectories $x^*(\cdot)$, but for all trajectories corresponding to control functions satisfying $J(x_0,u) \le V(x_0) + \delta$, i.e.\ for all \emph{near optimal trajectories}. For even more variants and a historic discussion of the turnpike property we refer to the recent survey \cite{FauG20}.

\subsection{Strict dissipativity and the turnpike property}

Under a boundedness assumption on the optimal value functions it is fairly easy to prove that strict dissipativity implies the turnpike property. The boundedness assumption demands that there are a $\KK_\infty$-function $\gamma$ and a constant $C>0$, such that the optimal value and the storage functions satisfy
\be |V_T(x)-V_T(x^e)| \le \gamma(\|x-x^e\|) + C \;\; \mbox{ and } \;\;  |\lambda(x)| \le \gamma(\|x-x^e\|) + C \label{eq:bounds}\ee
for all $x\in\X$ and all $T\ge 0$. The first condition can be ensured by a reachability condition: if the equilibrium $x^e$ can be reached from every $x\in\X$ with costs that are bounded in bounded subsets of $\X$, then the inequality for $|V_T(x)-V_T(x^e)|$ can be established for all $T>0$; for details see \cite[Section 6]{Grue13}.

We illustrate the derivation of the turnpike property under these boundedness assumptions in the discrete-time setting \LG{and refer to  \cite[Section 2]{FKJB17} for a continuous-time version of this result.} Using \eqref{eq:dtsdiss} and \eqref{eq:optsupply} and choosing $D>0$ such that $\lambda(x)\ge -D$ for all $x\in\X$ we obtain 
\beas J(x_0,u) & = &  \sum_{\LG{t=0}}^{T-1} \ell(x(t),u(t)) \\
& \ge & \sum_{\LG{t=0}}^{T-1} \alpha(\|x(t)-x^e\|) \LG{\,+\, T\ell(x^e,u^e)} - \lambda(x_0) + \lambda(x(T)) \\
& \ge & \sum_{\LG{t=0}}^{T-1} \alpha(\|x(t)-x^e\|) + T\ell(x^e,u^e) - \lambda(x_0) - D 
\eeas
Moreover, using the constant control $u\equiv u^e$ we obtain
\[ V_T(x^e) \le J_T(x^e,u) = \sum_{t=0}^{T-1} \ell(x^e,u^e) = T\ell(x^e,u^e).\]
Together this implies for $x=x^*$ with $x^*(0)=x_0$
\[ \gamma(\|x_0-x^e\|) + C \ge |V_T(x_0)-V_T(x^e)|  \ge \sum_{k=0}^{T-1} \alpha(\|x^*(t)-x^e\|) - \lambda(x_0) - D.\]
Given an arbitrary $\eps>0$, if $x^*$ spends too much time outside $B_\eps(x^e)$, then this inequality is violated. From this fact the existence of the constant $C$ in the turnpike property can be concluded. Here the inequalities in \eqref{eq:bounds} are needed in order to make sure that the constant $C$ in the turnpike property only depends on the size $K$ of the ball $B_K(x^e)$ containing $x_0$ and not on the individual initial value $x_0$. 

One easily checks that both examples \eqref{eq:brock} and \eqref{eq:invariance} are strictly dissipative, the former with the storage function $\lambda(x) = \alpha(x-x^e)/x^e$ and the latter with the storage function $\lambda(x) = -x^2/2$. This explains why we see the turnpike behavior in these examples. Note that in both examples boundedness of $\X$ is important to ensure that $\lambda$ is bounded from below on $\X$.

With a more involved proof, strict dissipativity with suitable bounds on the problem data can also be used in order to obtain an exponential turnpike property, see \cite{DGSW14}. Alternatively, the exponential turnpike property can be established using necessary optimality conditions (we refer to \cite{PorZ13,TreZ15,TrZZ18,GrSS20} for a selection of papers in this direction), but the strict dissipativity based approach has the advantage that it works for arbitrary nonlinearities, while the optimality condition-based approach is typically limited to linear systems or to small nonlinearities via linearization arguments.  \LG{The approach based on optimality conditions, in turn, has the advantage that it also provides a turnpike property for the adjoint states. For the dissipativity-based approach, such an implication is so far only known for the particular case of a so-called interval turpike property, in which the solutions are required to stay near the turnpike during a sufficiently large interval of time \cite{FGHS20}.}

The fact that strict dissipativity implies the turnpike property immediately raises the question how much stronger strict dissipativity is than the turnpike property. The answer lies in the observation that the inequality chain, above, can also be used for establishing the turnpike property for trajectories and controls that are not strictly optimal but only satisfy the inequality $J_T(x_0,u)\le V_T(x_0) + \delta$, i.e.\ they are near-optimal. In other words, we obtain a \emph{near-optimal turnpike property}. For systems that are locally controllable in a neighborhood of $x^e$ and controllable to this neighborhood from all $x\in\X$, it was shown in \cite{GruM16} that \LG{in discrete time} strict dissipativity is equivalent to the near-optimal turnpike property. \LG{In continuous time, to the best of the author's knowledge results in this generality are not yet known, but results under the stronger assumptions that the turnpike is exact or that the running cost $\ell$ is locally positive definite with respect to $(x^e,u^e)$ can be found in \cite[Section 4]{FKJB17}.}

Interestingly, the proof of the equivalence \LG{in discrete time} relies on a similar equivalence result for plain (i.e., nonstrict) dissipativity: \LG{In \cite{MuAA15} (for discrete-time systems) and in \cite{FKJB17} (for continuous-time systems)} it was shown that, again under a local controllability assumption, dissipativity is equivalent to the fact that the average cost for all admissible $u$ satisfies the inequality 
\be \LG{\liminf_{T\to\infty}} \frac1T J_T(x_0,u) \ge \ell(x^e,u^e), \label{eq:avbound}\ee
a property that is known under the name of \emph{optimal operation at steady state}.

Strict dissipativity also explains why optimal trajectories exhibiting the turnpike property typically look as in Figure \ref{fig:turnpike_examples} and hardly ever as in Figure \ref{fig:turnpike_illustration}. The reason is that a similar inequality as above shows that an excursion from the equilibrium $x^e$ followed by a return to $x^e$ causes costs that are significantly higher than staying in $x^e$. Under the assumption that the optimal values are continuous near $x^e$ (in the sense that the optimal value function near $x^e$ has about the same values as in $x^e$ --- a property that can again be rigorously established for locally controllable systems), this implies that such excursions will not happen in (near-)optimal trajectories.

In contrast to that, an excursion from the equilibrium $x^e$ can be cheaper then staying in $x^e$ if the solution does \emph{not} return to $x^e$ after the excursion. This is precisely the effect that creates the leaving arc of the optimal trajectories. As the gain that can be obtained from the leaving arc is bounded by the value of the storage function $\lambda(x(T))$ at the terminal state of the trajectory, it is important that $\lambda$ is bounded from below. Indeed, the example \eqref{eq:invariance} would cease to be strictly dissipative if we changed $\X=[-2,2]$ to $\X=\R$. In this case, it is easily seen that the optimal trajectories would tend to $\pm\infty$ (with optimal control $u^*\equiv 0$) instead of staying near $x^e=0$ most of the time, i.e.\ the turnpike property also ceases to exist.

\subsection{Related properties}

Strict dissipativity generalizes a lot of well-known properties for optimal control problems ensuring a certain asymptotic behavior for the optimal trajectories. To begin with, it is easily seen that strict dissipativity holds with storage function $\lambda\equiv 0$ for any cost function $\ell$ satisfying $\ell(x^e,u^e)=0$ and $\ell(x,u)\ge \alpha(\|x-x^e\|)$ for all $x\in\X$ and $u\in\U$. In case of a quadratic cost 
\[ \ell(x,u) = (x-x^e)^TQ(x-x^e) + (u-u^e)^TR(u-u^e)\]
this amounts to requiring that $Q$ is positive definite. For linear quadratic problems with dynamics $f(x,u)=Ax+Bu$, generalized quadratic cost 
\[ \ell(x,u) = x^TQx + u^TRu + q^Tx + r^Tu\]
with $Q=C^TC$, and no state constraints, i.e., $\X=\R^n$, strict dissipativity is equivalent to detectability of the pair $(A,C)$, i.e.\ to the fact that all unobservable eigenvalues \LG{$\mu$} of $A$ satisfy \LG{Re$\mu<0$} in continuous time or \LG{$|\mu|<1$} in discrete time. If $\X$ is bounded with $x^e$ in its interior, then strict dissipativity is equivalent to the fact that all unobservable eigenvalues \LG{$\mu$} of $A$ satisfy \LG{Re$\mu\ne0$} in continuous time or \LG{$|\mu|\ne 1$} in discrete time (see \cite{GruG18} and \cite{GruG21}  for proofs in discrete and continuous time, respectively). We note that the last criterion applies to example \eqref{eq:invariance}. \LG{While the results just cited apply to problems in which the $R$-matrix is positive definite, some of the implications were also be proved in the indefinite case, see \cite{BKAM18}. This paper also extends the analysis in \cite{GruG18} to quadratic costs including cross terms of the form $x^TSu$, by using a transformation that transforms a cost function with cross terms into a cost function without such terms. 
}

\LG{For nonlinear costs, it can be shown that strict dissipativity holds if $f$ is affine and $\ell$ is strictly convex, a criterion which applies to example \eqref{eq:brock}. In this case, the storage function can be chosen as a linear function $\lambda(x) = p^Tx$, where $p\in\R^n$ is the Lagrange multiplier of the optimization problem for determining the optimal equilibrium, i.e., in discrete time
\[ \min_{(x,u)\in \X\times U} \ell(x,u) \quad \mbox{s.t. } \;  f(x,u)=x.\] 
A proof can be found in \cite{DGSW14} and it is straightforward to carry over the statement and proof to continuous time. We remark that this fact was known and used before, e.g., in \cite{DiAR10}, but we were not able to find an earlier reference to a formal proof.  It should be noted that a linear function of the form $\lambda(x) = p^Tx$ is in general not bounded from below as required from a storage function. Hence, appropriate state constraints must be imposed in order to guarantee strict dissipativity. }

\LG{For fully nonlinear discrete-time systems}, as shown in \cite{HoeG19}, strict dissipativity also follows from nonlinear detectability notions, such as the one introduced in \cite{GMTT05} or input-output-to-state stability (IOSS) \cite{CaiT08}. \LG{We briefly explain the former connection. The nonlinear detectability condition from \cite{GMTT05} demands the existence of a nonnegative function $W:\X\to\R_{\ge 0}$ and of $\KK_\infty$ functions $\alpha_i$, $i=1,\ldots,3$ such that the inequalities 
\beas
W(x)& \le & \alpha_1(\|x-x^e\|)  \\
W(f(x,u)) - W(x) & \le & -\alpha_2(\|x-x^e\|) + \alpha_3(\ell(x,u))
\eeas
hold for all $x\in \X$ and $u\in\U$. Here one assumes $\ell(x,u)\ge 0$ for all $x\in\X$ and $u\in\U$ as well as $\ell(x^e,u^e)=0$, which is consistent with the fact that detectability applies to optimal control problems that are designed for stabilizing a given equilibrium, here $x^e$. Then, under a growth assumption on $\alpha_3$ (which holds without loss of generality if $\ell$ is bounded on $\X\times \U$), it was shown in \cite{HoeG19} that for an appropriately chosen $\rho\in\KK_\infty$ the function $\lambda = \rho\circ W$ is a storage function certifying strict dissipativity for the supply rate $s(x,u) = \ell(x,u)$. Conversely, if the problem is strictly dissipative and $\lambda(x^e)\le \lambda(x)$ for all $x\in\X$, then the problem is also detectable with $W(x)=\lambda(x)-\lambda(x^e)$. As mentioned, these results are in discrete time and we are not aware of the existence of continuous-time counterparts, although we conjecture that such result should be feasible. 
}

\subsection{\LG{Verifying (strict) dissipativity of optimal control problems}}

\LG{Since (strict) dissipativity of optimal control problems has many useful applications, it is, of course, desirable to be able to verify if an optimal control problem has this property. If none of the sufficient structural conditions discussed in the previous section can be used for this purpose, then other methods need to be applied. To this end, a couple of algorithmic approaches were developed in the last years.}

\LG{The conceptually simplest approach relies on the relation between strict dissipativity and the turnpike property. One can simply solve the optimal control problem numerically for different initial conditions and different time horizons and check whether the turnpike property holds, which is very easy to observe. This already yields a good indication whether strict dissipativity may hold, but, of course, it is only an indication and not a rigorous check.}

\LG{A rigorous check can be obtained by computing a storage function $\lambda$ numerically. Given the conceptual similarity between storage functions and Lyapunov functions, it seems plausible that most computational techniques for Lyapunov functions (for a survey see, e.g., \cite{GieH15}) can be modified for computing storage functions. A method for which this has already been done is the sum-of-squares approach, that was applied for the computation of storage functions in \cite{BKAM20} in continuous time and in \cite{PiAG19} in discrete time. The main bottleneck of this approach---and of most other numerical approaches for computing Lyapunov-like functions---is the fact that the computational effort of the methods grows rapidly with the space dimension of the problem. This is the so-called curse of dimensionality.}

\LG{For linear-quadratic problems this problem vanishes, because one can always resort to linear-quadratic storage functions, which can be efficiently computed using LMI techniques \cite{SchW00}. Recently, for this setting also a purely data-driven approach for verifying dissipativity was developed \cite{RBKA19}. For general nonlinear problems, a promising way appears to be the use of approximation architectures that can overcome the curse of dimensionality, such as deep neural networks. These architectures are provably able to compute approximate solutions of certain partial differential equations in high space dimensions \cite{DaLM20,HJKN20} and recently it was shown that the same is true for Lyapunov functions \cite{Grue20a}. Of course, this cannot work in general but only under suitable structural properties of the problem under consideration, but there is hope that these are related to well known system theoretic properties (such as small-gain conditions in case of \cite{Grue20a}), which are well understood. It is quite likely that we will soon see progress in this direction, particularly in connection with reinforcement learning approaches.}

\section{An application: Model Predictive Control}

Model Predictive Control (MPC), as described in the Sidebar entitled \nameref{sidebar:mpc}, is a highly popular control method, in which the computationally challenging solution of an infinite horizon optimal control problem is replaced by the successive solutions of optimal control problems on finite time horizons. As described in the sidebar, this induces an error that can be analyzed employing strict dissipativity of the optimal control problem.

\subsection{Stability analysis}

Strict dissipativity turns out to be useful for two important aspects of this analysis.  The first aspect is the stability analysis. In general, it is not clear, at all, that the MPC closed-loop solutions will exhibit stability-like behaviour. However, under the assumption that strict dissipativity holds, we can define an auxiliary optimal control problem by using the modified or rotated running or stage cost
\[ \tilde \ell(x,u)  := \ell(x,u) - \ell(x^e,u^e) - D\lambda(x)f(x,u)\]
in continuous time or 
\[ \tilde \ell(x,u)  := \ell(x,u) - \ell(x^e,u^e) + \lambda(x) - \lambda(f(x,u))\]
in discrete time. It is immediate from the strict dissipativity inequalities \eqref{eq:ctdinfs} or \eqref{eq:dtdones} that this modified stage cost is positive definite at $x^e$, i.e., it satisfies $\tilde\ell(x^e,u^e)=0$ and $\tilde \ell(x,u) \ge \alpha(\|x-x^e\|)$. Thus, if we consider the optimal control problems \eqref{eq:ctoc} or \eqref{eq:dtoc} with $\tilde\ell$ in place of $\ell$, it forces the optimal solutions to approach $x^e$. Denote the corresponding functional and optimal value function by $\widetilde J_T$ and $\widetilde V_T$, respectively, and the optimal solution by $\tilde x^*(\cdot)$. Then by either adding a terminal cost $\widetilde F$ with suitable properties to $\widetilde J_T$ (cf., e.g., \cite{CheA98,MaRRS00}, \cite[Chapter 2]{RaMD17}, or \cite[Chapter 5]{GruP17}) or by imposing conditions that ensure that $\widetilde V_T$ is close to $\widetilde V_\infty$ for sufficiently large $T$ (cf., e.g., \cite{GMTT05} or \cite[Chapter 6]{GruP17}), one can prove that $\widetilde V_T$ is a Lyapunov function for the MPC closed loop. The trouble is, though, that this is only true if $\tilde\ell$ is used as running cost or stage cost in the MPC scheme. This, however, is often not practical because the storage function may not be known and difficult to compute. In this case, one would  like to keep the original cost $\ell$ in the optimal control problem. 

In case the running or stage cost \LG{$\ell$} together with the terminal cost \LG{$F$} is used, the trick is to define \LG{$\tilde\ell$ as above and $\widetilde F := F + \lambda$. Then the optimal control problems using $\ell$ and $F$} produces the same optimal solutions as the one using $\tilde \ell$ and $\widetilde F$, \LG{and the analysis described in the previous paragraph can be used, without the need to actually compute $\tilde \ell$ and $\widetilde F$. This is because} the conditions on the terminal cost needed to make $\widetilde V_T$ a Lyapunov functions can be formulated directly for $F$ without having to use (or even to know) $\lambda$. The resulting required inequality for $F$ is of the form 
\be DF(x)f(x,u) \le - \ell(x,u) + \ell(x^e,u^e) \label{eq:Fct} \ee
in continuous time and 
\be F(f(x,u)) \le F(x) - \ell(x,u) + \ell(x^e,u^e) \label{eq:Fdt} \ee
in discrete time, which must hold for all $x$ in the terminal constraint set $\X_0$ and a control value $u$ (depending on $x$), which is such that the solution does not leave $\X_0$. 
Since the optimal control problems with $\tilde\ell$ and $\widetilde F$ on the one side and with $\ell$ and $F$ on the other side produce the same optimal solutions, the MPC closed-loops resulting from the two problems also coincide and $\widetilde V_T$ can again be used as a Lyapunov function to conclude asymptotic stability. This trick was first proposed in \cite{DiAR10} and then refined in \cite{AnAR12}.

\begin{figure}[htb]
\includegraphics[width=7cm]{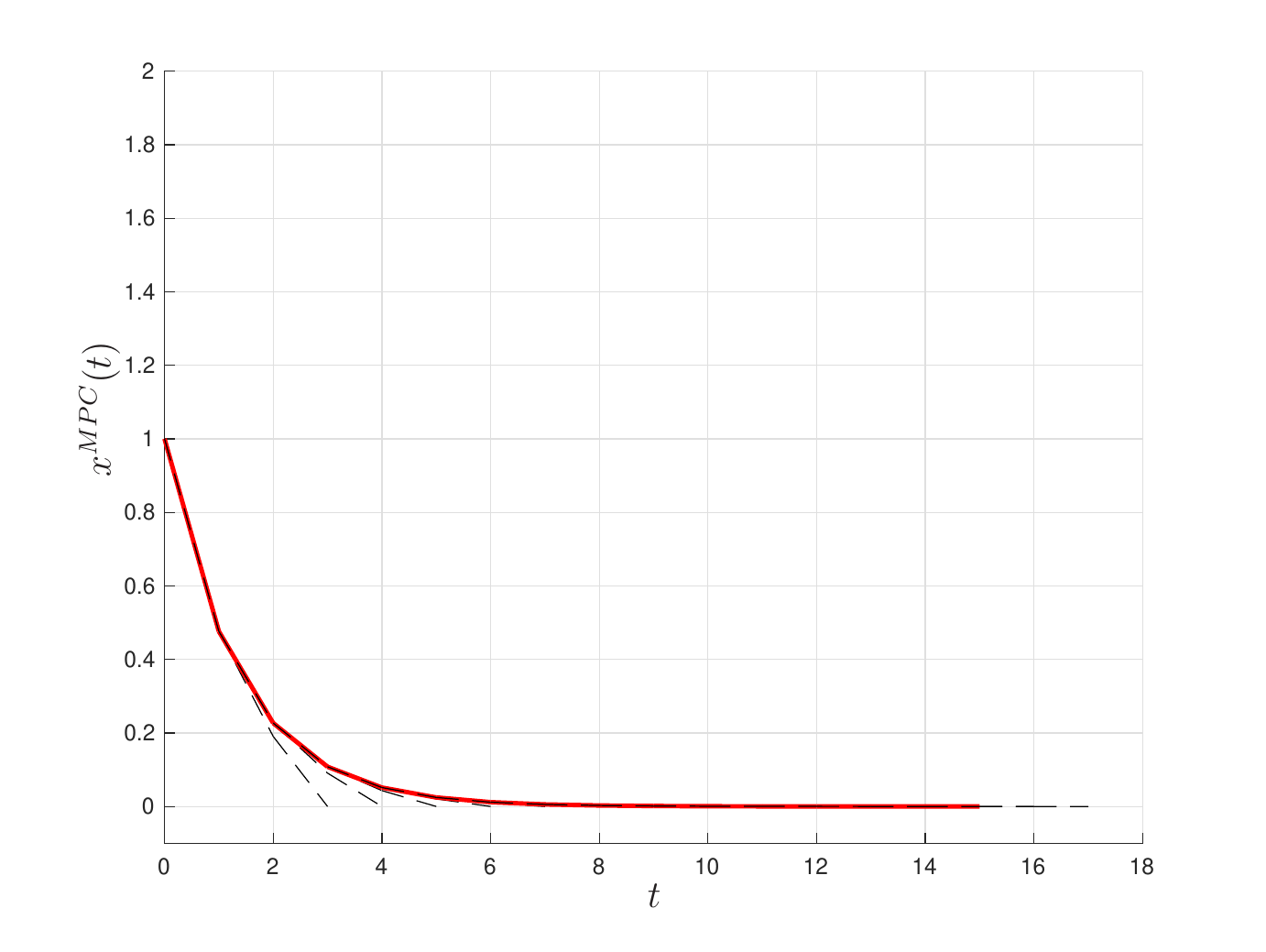}
\includegraphics[width=7cm]{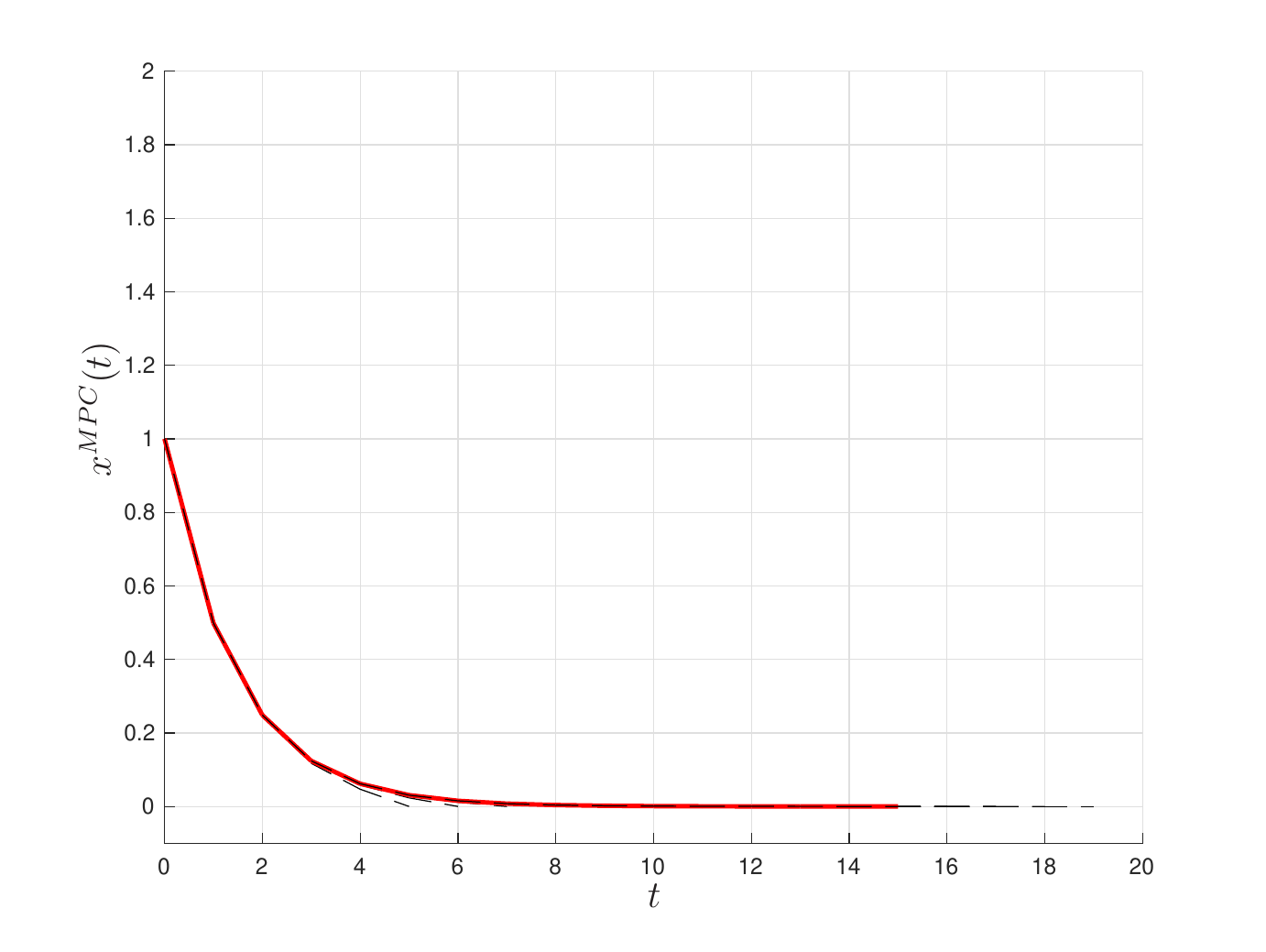}
\caption{MPC closed-loop trajectories (red solid) and predictions (black dashed) for example \eqref{eq:invariance}, using the terminal constraint set $\X_0=\{0\}$ and the terminal cost $F\equiv 0$ and horizon $T=3$ (left) and $T=5$ (right). Note that the terminal constraint $x(T)=0$ forces all predictions to end in $x^e=0$. \label{fig:inv_tc}}
\end{figure}

Figure \ref{fig:inv_tc} \LG{shows} the MPC closed-loop trajectories and the corresponding predictions for example \eqref{eq:invariance}, using the terminal constraint set $\X_0=\{0\}$ and the terminal cost $F\equiv 0$. 

In many cases it may be difficult to find a terminal cost $F$ that meets the required condition \eqref{eq:Fct} or \eqref{eq:Fdt}. While the trivial choice $\X_0=\{x^e\}$ and $F\equiv 0$ always satisfies \eqref{eq:Fct} or \eqref{eq:Fdt}, this choice of $\X_0$ may cause problems in the numerical optimization and result in a small set of feasible states for \eqref{eq:ctoc} or \eqref{eq:dtoc} when the terminal condition $x(T)\in \X_0$ is added. It may thus be attractive to drop the terminal constraint $\X_0$ and costs $F$ \LG{and $\widetilde F$}. In this case, however, \LG{the trick enforcing identical solutions of the optimal control problems with $\ell$ and with $\tilde \ell$ by setting $\widetilde F := F + \lambda$ is no longer applicable}. In fact, without terminal costs the optimal trajectories and controls of the problem with $\ell$ and with $\tilde \ell$ do not coincide anymore. However, they may still coincide approximately in a suitable sense.

In order to establish this property, we need to assume that the storage function $\lambda$ and the optimal value functions $V_T$ and $\widetilde V_T$ are continuous in $x^e$, uniformly in $T$ (the former is not very restrictive since often $\lambda$ is a polynomial and the latter can again be ensured by local controllability, see \cite[Section 6]{Grue13}). Next we use that both the problem with cost $\ell$ and the problem with cost $\tilde \ell$ are strictly dissipative and thus exhibit the turnpike property. From this we can conclude that for any two optimal trajectories $x^*(\cdot)$ and $\tilde x^*(\cdot)$ starting in the same initial value $x_0$, there will be a time $\tau$ such that they satisfy $\tilde x^*(\tau)\approx x^e$ and $x^*(\tau)\approx x^e$, where the error hidden in ``$\approx$'' tends to $0$ as the horizon $T$ increases. Together with the continuity assumption on the optimal value function this implies that the cost of the two trajectories on the time interval $[0,\tau]$ is almost identical. This, finally, can be used to conclude that $\widetilde V_T$ is an approximate Lyapunov function, from which practical asymptotic stability, i.e., asymptotic stability of a neighborhood of $x^e$, which shrinks down to $\{x^e\}$ when $T$ increases, can finally be concluded. The details of this reasoning were originally derived in the papers \cite{Grue13,GruS14} \LG{in discrete time} and in \cite{FauB15} \LG{in continuous time}. A concise presentation can be found in \cite[Section 8.6]{GruP17} or in \cite[Section 4]{FaGM18}. 

\begin{figure}[htb]
\includegraphics[width=7cm]{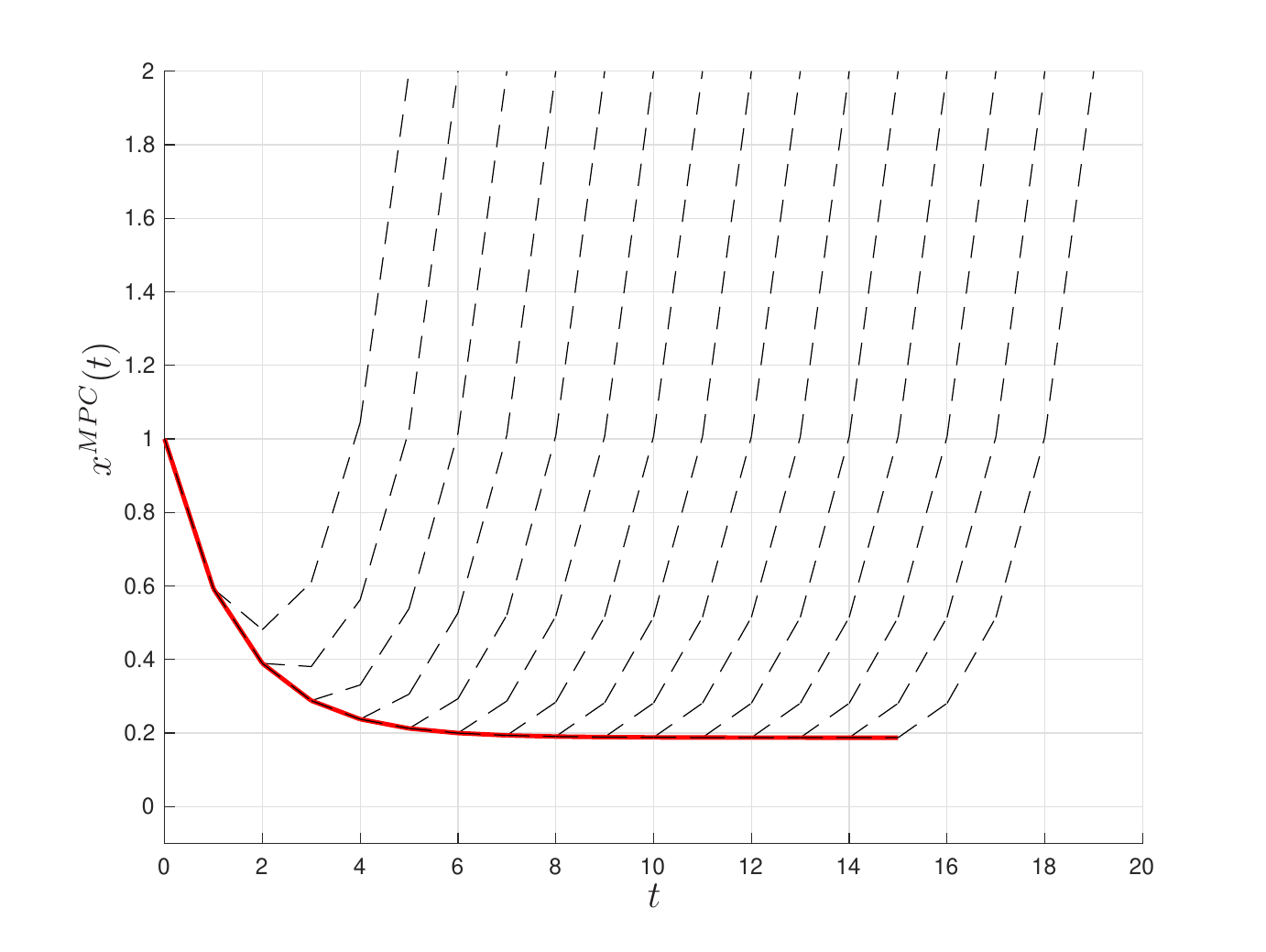}
\includegraphics[width=7cm]{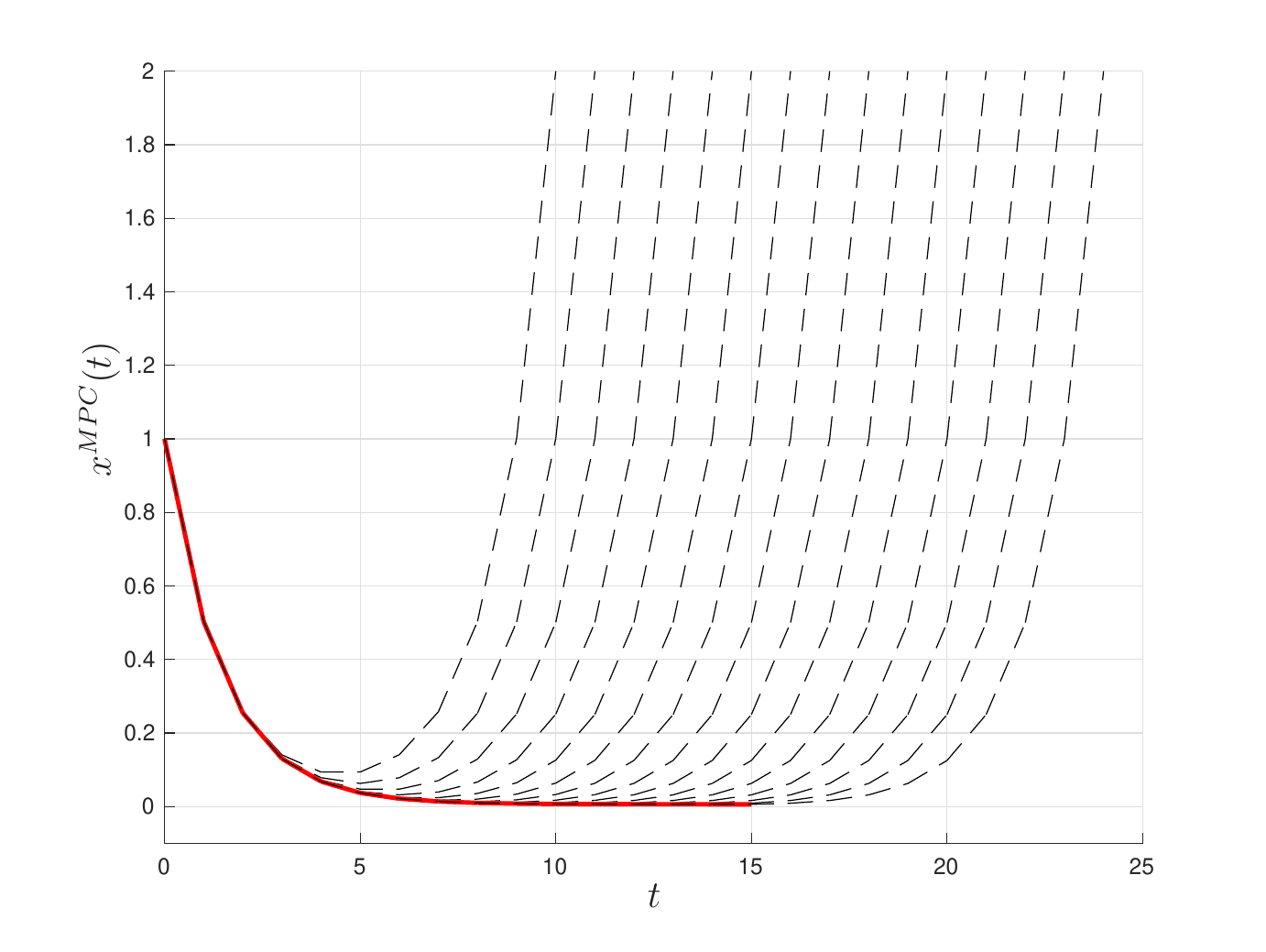}
\caption{MPC closed-loop trajectories (red solid) and predictions (black dashed) for example \eqref{eq:invariance} without terminal conditions with horizon $T=5$ (left) and $T=10$ (right). Without any terminal constraints all predictions end in $x=2$, cf.\ also Figure \ref{fig:turnpike_examples}(right). \label{fig:inv_notc}}
\end{figure}

Figure \ref{fig:inv_notc} show the MPC closed-loop trajectories and the corresponding predictions for example \eqref{eq:invariance} without any terminal conditions. In comparison to Figure \ref{fig:inv_tc} the merely practical asymptotic stability of the equilibrium $x^e=0$ is clearly visible, since the red closed-loop solution does not converge to the equilibrium $x^e=0$ but only to a neighborhood of this equilibrium, which becomes smaller as $T$ increases.

We note that under stronger assumptions than strict dissipativity stronger stability statements can be made for MPC without terminal conditions. For instance, under \LG{the nonlinear detectability condition discussed in the previous section} and suitable bounds on the optimal value function, it was shown in \cite{GMTT05} that true (as opposed to merely practical) asymptotic stability of $x^e$ for the MPC closed loop can be expected for sufficiently large $T$. For positive definite cost, different ways of estimating the length of the horizon $T$ needed for obtaining asymptotic stability were proposed in \cite{TuMT06,Grue09,GPSW10,RebA11} (see also \cite[Chapter 6]{GruP17}). 
\LG{For linear quadratic problems, it was shown in \cite{ZaGD14} that for indefinite problems that yield asymptotically stable optimal solutions it is always possible to find an equivalent LQ problem with positive definite costs, to which these results are applicable.}
\LG{Finally, under suitable conditions true asymptotic stability can also be enforced using linear terminal costs that are much simpler than those in \eqref{eq:Fct} or  \eqref{eq:Fdt} and that are derived from the adjoints of the optimal control problem, see \cite{ZanF18}.}

\subsection{Performance analysis}

Since MPC relies on the solution of optimal control problems, it seems reasonable to expect that the MPC closed loop also enjoys certain optimality properties. \LG{This is the second aspect of the analysis of MPC schemes where dissipativity is helpful.}  To this end, we denote the MPC closed-loop trajectory by $x_{MPC}(t)$ and the corresponding control by $u_{MPC}(t)$. Then we can define different measures for the closed loop cost: 
The infinite horizon performance
\[ J^{MPC}_\infty(x_0) := \sum_{t=0}^\infty \ell(x_{MPC}(t),u_{MPC}(t)) \] 
would be the ``natural'' measure if we consider MPC as an approximation to an infinite horizon problem. 
However, as the infinite sum may not converge, we also look at other measures. We also consider the finite horizon closed-loop performance 
\begin{equation} J^{MPC}_S(x_0) := \sum_{t=0}^{S-1} \ell(x_{MPC}(t),u_{MPC}(t)) \label{eq:JmuK}\end{equation} 
and the averaged infinite horizon performance 
\[ \overline J^{MPC}_\infty(x_0) := \limsup_{S\to\infty}\frac 1S J^{MPC}_S(x_0).\] 
These last two performance measures complement each other, as the first measures the performance on finite intervals $[0,S]$ while the second measures the performance in the limit for $S\to\infty$.

The derivation of estimates for these quantities heavily relies on strict dissipativity and the turnpike property, which are exploited for MPC both with and without terminal conditions, and on the stability properties described in the previous section. The key idea is to use the similarity of the initial pieces of optimal trajectories until they reach the optimal equilibrium in order to derive approximate versions of the dynamic programming principle. From these, estimates on the above quantities can be derived by induction over $t$. The following \LG{discrete-time} estimates were originally developed in the papers \cite{Grue13,GruS14,GruP15a}. Concise presentations can be found in \cite{Grue16} or \cite[Chapter 8]{GruP17}. \LG{In continuous-time, to the best of our knowledge so far only results for the averaged performance are available. These can be found in \cite{AlAJ16}.}

For MPC with terminal conditions satisfying \eqref{eq:Fdt}, the identity
\be \overline J^{MPC}_\infty(x_0) = \ell(x^e,u^e) \label{eq:av_tc}\ee
holds for all $N\in\N$, and because of \eqref{eq:avbound} this is the best possible value the average performance functional can attain. In case $V_\infty$ assumes finite values, there exists a function $\delta_1:\N\to[0,\infty)$ with $\delta_1(T)\to 0$ as $T\to\infty$ such that the inequality
\be J^{MPC}_\infty(x_0) \le V_\infty(x_0) + \delta_1(T) \label{eq:inf_tc}\ee
holds. In case $V_\infty$ does not assume finite values, for each $S\in\N$ we can obtain the estimate
\be J^{MPC}_S(x_0) \le \inf_{u\in \widetilde \U^S} J_S(x_0,u) + \delta_1(T) + \delta_2(S) \label{eq:transient_tc}\ee
where $\delta_2$ is the same type of function as $\delta_1$. Here $\widetilde\U^S$ denotes the set of admisible controls for which $\|x_u(S,x_0)-x^e\| \le \|x^{MPC}(S,x_0)-x^e\|$ holds. As $x^{MPC}(S,x_0) \to x^e$ holds for $S\to\infty$, for large $S$ the quantity $\inf_{u\in \widetilde \U^S} J_S(x_0,u)$ measures the optimal transient cost for trajectory going from $x_0$ to a small neighborhood of $x^e$. Thus, the estimates show that MPC produces trajectories that produce optimal averaged cost and approximately optimal transient cost. 

\begin{figure}[htb]
\includegraphics[width=7cm]{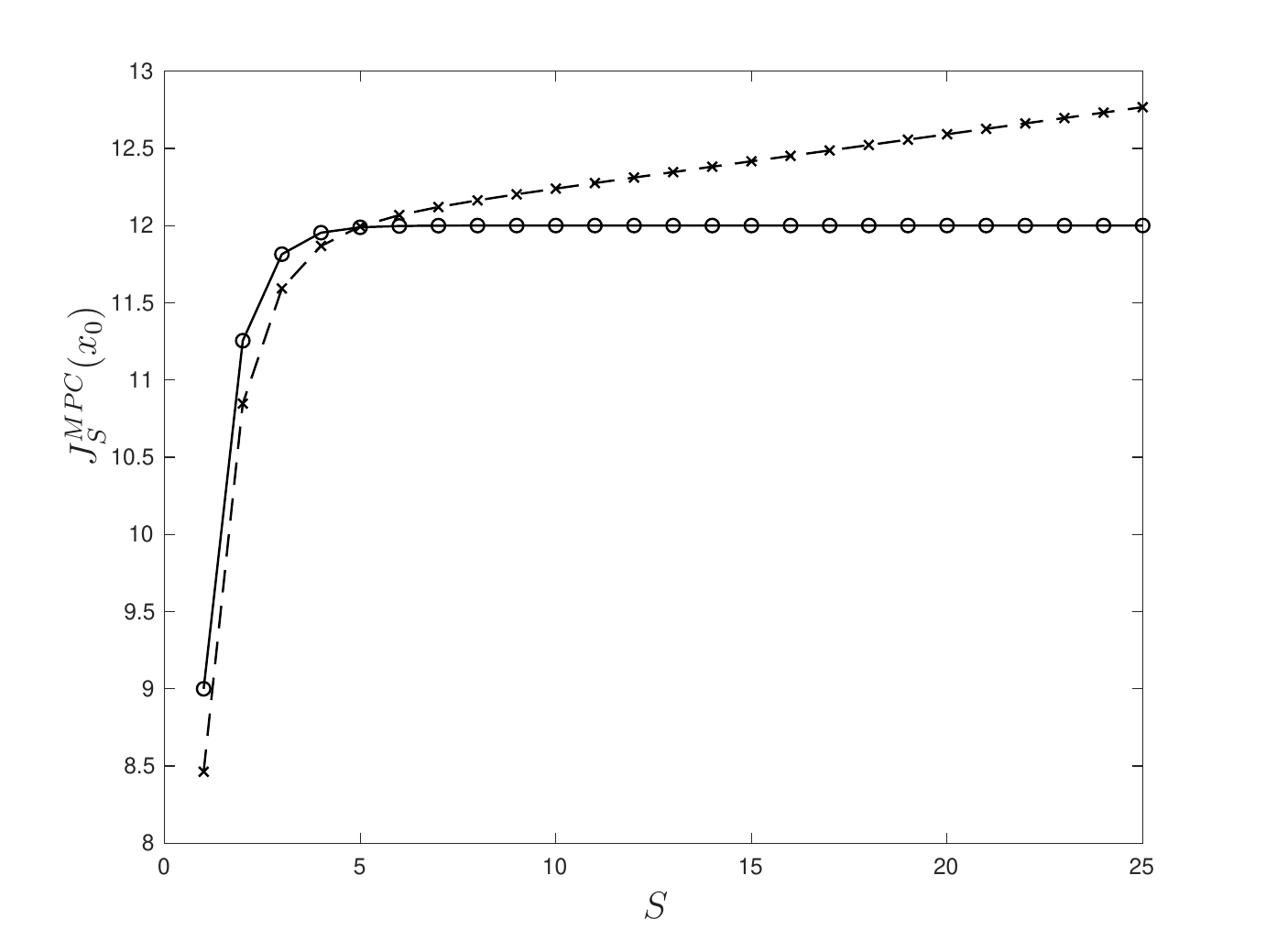}
\includegraphics[width=7cm]{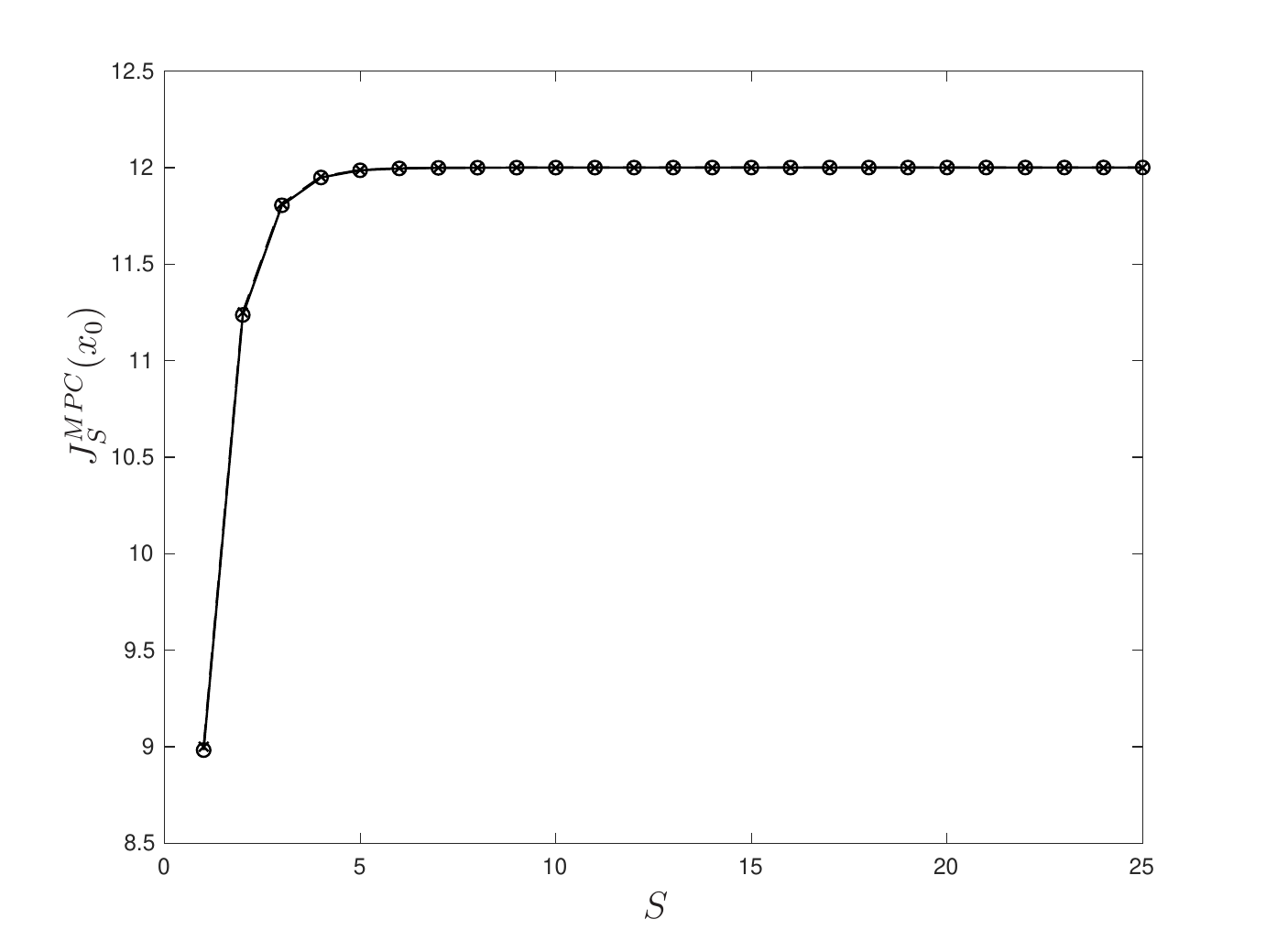}
\caption{MPC closed-loop cost $J_S^{MPC}(x_0)$ for example \eqref{eq:invariance} with $x_0=2$ for varying $S$. The solid line with circles shows the values with terminal conditions $\X_0=\{0\}$ and $F\equiv 0$, the dashed line with crosses the values without terminal conditions, with horizon $T=5$ (left) and $T=10$ (right).\label{fig:J_compare}}
\end{figure}

Without terminal conditions the results become somewhat weaker. Particularly, we can in general no longer ensure that $J^{MPC}_\infty(x_0)$ is finite, even if $V_\infty(x_0)$ is finite. However, we can still establish counterparts of \eqref{eq:av_tc} and \eqref{eq:transient_tc}, namely
\be \overline J^{MPC}_\infty(x_0) \le \ell(x^e,u^e) + \delta_1(T) \label{eq:av_noatc}\ee
and
\be J^{MPC}_S(x_0) \le \inf_{u\in \widetilde \U^S} J_S(x_0,u) + S\delta_1(T) + \delta_2(S), \label{eq:transient_notc}\ee
with $\delta_1$ and $\delta_2$ of the same type as above. We note that the fact that the error term $\delta_1(T)$ in \eqref{eq:transient_tc} increases to $S\delta_1(T)$ in \eqref{eq:transient_notc} is not an effect of an insufficiently precise analysis but actually a natural consequence of the mere practical asymptotic stability of $x^e$ without terminal conditions: If the closed-loop solution does not converge to $x^e$, the stage cost will typically not converge to $0$ but to a nonzero residual value, which keeps accumulating over the time $S$. In order to illustrate this effect, Figure \ref{fig:J_compare} shows the values of $J^{MPC}_S(x_0)$ for example \eqref{eq:invariance} with terminal conditions ($\circ$) and without terminal conditions ($\times$) for different $S$. In the left figure with horizon $T=5$, the effect of the additional factor $S$ in front of $\delta_1(T)$ in \eqref{eq:transient_notc} is clearly visible. In the right figure with horizon $T=10$, the error term $\delta_1(T)$ in \eqref{eq:transient_notc} is already so small that the effect of the factor $S$ is not visible anymore. For this horizon, MPC with and without terminal conditions yield almost the same performance. 

Just as for stability, under stronger assumptions stronger statements can be made. For instance, for positive definite cost $\ell$ finiteness of $J^{MPC}_\infty(x_0)$ and an estimate of the form \eqref{eq:inf_tc} can be obtained also for MPC without terminal conditions, see \cite{GPSW10,RebA11}.

\section{Summary, extensions and outlook}

Dissipativity and strict dissipativity are important systems theoretic properties with multiple applications. In this paper we have shown that they naturally link to optimal control, a fact that is already prominently present in Jan C.\ Willems' earliest publications on the subject. Here we provided a survey on recent results in this direction, which establish a close link between strict dissipativity and the turnpike property and showed how these concepts can be used for analyzing stability and performance of MPC schemes. 

The present results have already been extended into various directions. In particular, dissipativity has been extended to optimal control problems that do not exhibit an optimal equilibrium but an optimal periodic orbit \cite{ZaGD13,GruM16a,ZaGD17,KoMA18} or general time-varying optimal trajectories \cite{GrPS18,GruP19}. In this context, the concept of overtaking optimality can be used in order to define a meaning for optimality also in the case that the infinite horizon optimal value function $V_\infty$ is not finite. 

\LG{In the field of dissipativity and optimal control there are still a variety of open questions and we end by mentioning some of them. For instance, as pointed out throughout the paper, several important results are so far only available in discrete time and it would be very interesting to see if and how they can be carried over to continuous time.}

\LG{Another} major open question is the relation between strict dissipativity and detectability-like notions for infinite-dimensional systems. Even for linear-quadratic optimal control problems this relation is not yet fully understood. 

\LG{Finally, it is known that the turnpike property also appears in various forms in stochastic optimal control problems, see, e.g., the recent numerical study \cite{OBGF20}. While some results on turnpikes in the stochastic case are available, cf.\ \cite{Mari89,SSZJ92,KolY12}, the connection to dissipativity concepts is so far largely unexplored.}

\bibliography{../../bib/lars_consolidated}

\begin{thebibliography}{10}

\bibitem{AlAJ16}
{\sc A.~Alessandretti, A.~P. Aguiar, and C.~N. Jones}, {\em On convergence and
  performance certification of a continuous-time economic model predictive
  control scheme with time-varying performance index}, Automatica, 68 (2016),
  pp.~305--313.

\bibitem{AnAR12}
{\sc D.~Angeli, R.~Amrit, and J.~B. Rawlings}, {\em On average performance and
  stability of economic model predictive control}, IEEE Trans. Autom. Control,
  57 (2012), pp.~1615--1626.

\bibitem{BKAM18}
{\sc J.~Berberich, J.~K\"{o}hler, F.~Allg\"{o}wer, and M.~A. M\"{u}ller}, {\em
  Indefinite linear quadratic optimal control: strict dissipativity and
  turnpike properties}, IEEE Control Syst. Lett., 2 (2018), pp.~399--404.

\bibitem{BKAM20}
{\sc J.~Berberich, J.~K\"{o}hler, F.~Allg\"{o}wer, and M.~A. M\"{u}ller}, {\em
  Dissipativity properties in constrained optimal control: a computational
  approach}, Automatica, 114 (2020), pp.~108840, 9.

\bibitem{BreP20}
{\sc T.~Breiten and L.~Pfeiffer}, {\em On the turnpike property and the
  receding-horizon method for linear-quadratic optimal control problems}, SIAM
  J. Control Optim., 58 (2020), pp.~1077--1102.

\bibitem{BroM72}
{\sc W.~A. Brock and L.~Mirman}, {\em Optimal economic growth and uncertainty:
  the discounted case}, J. Econ. Theory, 4 (1972), pp.~479--513.

\bibitem{ByrL94}
{\sc C.~I. Byrnes and W.~Lin}, {\em Losslessness, feedback equivalence, and the
  global stabilization of discrete-time nonlinear systems}, IEEE Trans.
  Automat. Control, 39 (1994), pp.~83--98.

\bibitem{CaiT08}
{\sc C.~Cai and A.~R. Teel}, {\em Input-output-to-state stability for
  discrete-time systems}, Automatica, 44 (2008), pp.~326--336.

\bibitem{CaHL91}
{\sc D.~A. Carlson, A.~B. Haurie, and A.~Leizarowitz}, {\em Infinite horizon
  optimal control --- Deterministic and Stochastic Systems}, Springer-Verlag,
  Berlin, second~ed., 1991.

\bibitem{CheA98}
{\sc H.~Chen and F.~Allg{\"o}wer}, {\em A quasi-infinite horizon nonlinear
  model predictive control scheme with guaranteed stability}, Automatica, 34
  (1998), pp.~1205--1217.

\bibitem{DGSW14}
{\sc T.~Damm, L.~Gr\"une, M.~Stieler, and K.~Worthmann}, {\em An exponential
  turnpike theorem for dissipative discrete time optimal control problems},
  SIAM J. Control Optim., 52 (2014), pp.~1935--1957.

\bibitem{DaLM20}
{\sc J.~Darbon, G.~P. Langlois, and T.~Meng}, {\em Overcoming the curse of
  dimensionality for some {H}amilton-{J}acobi partial differential equations
  via neural network architectures}, Res. Math. Sci., 7 (2020), pp.~Paper No.
  20, 50.

\bibitem{DiAR10}
{\sc M.~Diehl, R.~Amrit, and J.~B. Rawlings}, {\em A {L}yapunov function for
  economic optimizing model predictive control}, IEEE Trans. Autom. Control, 56
  (2011), pp.~703--707.

\bibitem{DoSS58}
{\sc R.~Dorfman, P.~A. Samuelson, and R.~M. Solow}, {\em Linear programming and
  economic analysis}, A Rand Corporation Research Study, McGraw-Hill, New
  York-Toronto-London, 1958.

\bibitem{FauB15}
{\sc T.~Faulwasser and D.~Bonvin}, {\em On the design of economic {N}{M}{P}{C}
  based on approximate turnpike properties}, in Proceedings of the 54th IEEE
  Conference on Decision and Control --- CDC 2015, 2015, pp.~4964--4970.

\bibitem{FauG20}
{\sc T.~Faulwasser and L.~Gr\"une}, {\em Turnpike properties in optimal
  control: an overview of discrete-time and continuous-time results}.
\newblock To appear in the Handbook of Numerical Analysis, 2021.
\newblock Preprint available from https://arxiv.org/pdf/2011.13670.pdf.

\bibitem{FGHS20}
{\sc T.~Faulwasser, L.~Gr{\"u}ne, J.-P. Humaloja, and M.~Schaller}, {\em The
  interval turnpike property for adjoints}.
\newblock arXiv Preprint 2005.12120, 2020.

\bibitem{FaGM18}
{\sc T.~Faulwasser, L.~Gr\"une, and M.~A. M\"uller}, {\em Economic nonlinear
  model predictive control}, Foundations and
  Trends\textsuperscript{\textregistered} in Systems and Control, 5 (2018),
  pp.~1--98.

\bibitem{FauK20}
{\sc T.~Faulwasser and C.~M. Kellett}, {\em On continuous-time infinite horizon
  optimal control---dissipativity, stability and transversality}.
\newblock arXiv Preprint 2001.09601, 2020.

\bibitem{FKJB17}
{\sc T.~Faulwasser, M.~Korda, C.~N. Jones, and D.~Bonvin}, {\em On turnpike and
  dissipativity properties of continuous-time optimal control problems},
  Automatica, 81 (2017), pp.~297--304.

\bibitem{FPHG15}
{\sc M.~G. Forbes, R.~S. Patwardhan, H.~Hamadah, and R.~B. Gopaluni}, {\em
  Model predictive control in industry: Challenges and opportunities}, in
  Proceedings of the 9th IFAC Symposium on Advanced Control of Chemical
  Processes --- ADCHEM 2015, vol.~48 of IFAC-PapersOnLine, Whistler, Canada,
  2015, pp.~531--538.

\bibitem{FreK96}
{\sc R.~A. Freeman and P.~V. Kokotovi\'{c}}, {\em Robust nonlinear control
  design}, Systems \& Control: Foundations \& Applications, Birkh\"{a}user
  Boston, Inc., Boston, MA, 1996.
\newblock State-space and Lyapunov techniques.

\bibitem{GieH15}
{\sc P.~Giesl and S.~Hafstein}, {\em Review on computational methods for
  {L}yapunov functions}, Discrete Contin. Dyn. Syst. Ser. B, 20 (2015),
  pp.~2291--2331.

\bibitem{GMTT05}
{\sc G.~Grimm, M.~J. Messina, S.~E. Tuna, and A.~R. Teel}, {\em Model
  predictive control: for want of a local control {L}yapunov function, all is
  not lost}, IEEE Trans. Automat. Control, 50 (2005), pp.~546--558.

\bibitem{Grue09}
{\sc L.~Gr\"une}, {\em Analysis and design of unconstrained nonlinear {M}{P}{C}
  schemes for finite and infinite dimensional systems}, SIAM J. Control Optim.,
  48 (2009), pp.~1206--1228.

\bibitem{Grue13}
{\sc L.~Gr{\"u}ne}, {\em Economic receding horizon control without terminal
  constraints}, Automatica, 49 (2013), pp.~725--734.

\bibitem{Grue16}
{\sc L.~Gr\"une}, {\em Approximation properties of receding horizon optimal
  control}, Jahresber. DMV, 118 (2016), pp.~3--37.

\bibitem{Grue20a}
{\sc L.~Gr\"une}, {\em Computing {L}yapunov functions using deep neural
  networks}, J. Comput. Dyn., Online First (2020), p.~22 pages.

\bibitem{GruG18}
{\sc L.~Gr\"une and R.~Guglielmi}, {\em Turnpike properties and strict
  dissipativity for discrete time linear quadratic optimal control problems},
  SIAM J. Cont. Optim., 56 (2018), pp.~1282--1302.

\bibitem{GruG21}
{\sc L.~Gr\"une and R.~Guglielmi}, {\em On the relation between turnpike
  properties and dissipativity for continuous time linear quadratic optimal
  control problems}, Math. Control Rel. Fields, 11 (2021), pp.~169--188.

\bibitem{GrKW17}
{\sc L.~Gr\"une, C.~M. Kellett, and S.~R. Weller}, {\em On the relation between
  turnpike properties for finite and infinite horizon optimal control
  problems}, J. Optim. Theory Appl., 173 (2017), pp.~727--745.

\bibitem{GruM16}
{\sc L.~Gr\"une and M.~A. M\"uller}, {\em On the relation between strict
  dissipativity and the turnpike property}, Syst. Contr. Lett., 90 (2016),
  pp.~45--53.

\bibitem{GruP15a}
{\sc L.~Gr\"une and A.~Panin}, {\em On non-averaged performance of economic
  {M}{P}{C} with terminal conditions}, in Proceedings of the 54th IEEE
  Conference on Decision and Control --- CDC 2015, Osaka, Japan, 2015,
  pp.~4332--4337.

\bibitem{GruP17}
{\sc L.~Gr{\"u}ne and J.~Pannek}, {\em Nonlinear Model Predictive Control.
  Theory and Algorithms}, Springer-Verlag, London, 2nd~ed., 2017.

\bibitem{GPSW10}
{\sc L.~Gr\"une, J.~Pannek, M.~Seehafer, and K.~Worthmann}, {\em Analysis of
  unconstrained nonlinear {M}{P}{C} schemes with time varying control horizon},
  SIAM J. Control Optim., 48 (2010), pp.~4938--4962.

\bibitem{GruP19}
{\sc L.~Gr{\"u}ne and S.~Pirkelmann}, {\em Economic model predictive control
  for time-varying system: Performance and stability results}, Opt. Control
  Appl. Meth., 41 (2019), pp.~42--64.
\newblock Special Issue: MPC for Energy Systems: Economic and Distributed
  Approaches.

\bibitem{GrPS18}
{\sc L.~Gr\"{u}ne, S.~Pirkelmann, and M.~Stieler}, {\em Strict dissipativity
  implies turnpike behavior for time-varying discrete time optimal control
  problems}, in Control systems and mathematical methods in economics: Essays
  in Honor of Vladimir M. Veliov, G.~Feichtinger, R.~M. Kovacevic, and
  G.~Tragler, eds., vol.~687 of Lecture Notes in Econom. and Math. Systems,
  Springer, Cham, 2018, pp.~195--218.

\bibitem{GrSS19}
{\sc L.~Gr\"une, M.~Schaller, and A.~Schiela}, {\em Sensitivity analysis of
  optimal control for a class of parabolic {P}{D}{E}s motivated by model
  predictive control}, SIAM J. Control Optim., 57 (2019), pp.~2753--2774.

\bibitem{GrSS20}
{\sc L.~Gr\"une, M.~Schaller, and A.~Schiela}, {\em Exponential sensitivity and
  turnpike analysis for linear quadratic optimal control of general evolution
  equations}, J. Differ. Equ., 268 (2020), pp.~7311--7341.

\bibitem{GruS14}
{\sc L.~Gr\"une and M.~Stieler}, {\em Asymptotic stability and transient
  optimality of economic {MPC} without terminal conditions}, J. Proc. Control,
  24 (2014), pp.~1187--1196.

\bibitem{GuTZ16}
{\sc M.~Gugat, E.~Tr\'elat, and E.~Zuazua}, {\em Optimal {N}eumann control for
  the 1{D} wave equation: finite horizon, infinite horizon, boundary tracking
  terms and the turnpike property}, Systems Control Lett., 90 (2016),
  pp.~61--70.

\bibitem{HoeG19}
{\sc M.~H\"oger and L.~Gr\"une}, {\em On the relation between detectability and
  strict dissipativity for nonlinear discrete time systems}, IEEE Control Syst.
  Lett., 3 (2019), pp.~458--462.

\bibitem{HJKN20}
{\sc M.~Hutzenthaler, A.~Jentzen, T.~Kruse, and T.~A. Nguyen}, {\em A proof
  that rectified deep neural networks overcome the curse of dimensionality in
  the numerical approximation of semilinear heat equations}, SN Partial Differ.
  Equ. Appl., 10 (2020), p.~34.

\bibitem{KoMA18}
{\sc J.~K\"ohler, M.~A. M\"uller, and F.~Allg\"ower}, {\em On periodic
  dissipativity notions in economic model predictive control}, IEEE Control
  Syst. Lett., 2 (2018), pp.~501--506.

\bibitem{KolY12}
{\sc V.~Kolokoltsov and W.~Yang}, {\em Turnpike theorems for {M}arkov games},
  Dyn. Games Appl., 2 (2012), pp.~294--312.

\bibitem{Mari89}
{\sc R.~Marimon}, {\em Stochastic turnpike property and stationary
  equilibrium}, J. Econom. Theory, 47 (1989), pp.~282--306.

\bibitem{MaRRS00}
{\sc D.~Q. Mayne, J.~B. Rawlings, C.~V. Rao, and P.~O.~M. Scokaert}, {\em
  Constrained model predictive control: stability and optimality}, Automatica,
  36 (2000), pp.~789--814.

\bibitem{MoyA73}
{\sc P.~J. Moylan and B.~D.~O. Anderson}, {\em Nonlinear regulator theory and
  an inverse optimal control problem}, IEEE Trans. Automatic Control, AC-18
  (1973), pp.~460--465.

\bibitem{MuAA15}
{\sc M.~A. M\"uller, D.~Angeli, and F.~Allg\"ower}, {\em On necessity and
  robustness of dissipativity in economic model predictive control}, IEEE
  Trans. Autom. Control, 60 (2015), pp.~1671--1676.

\bibitem{GruM16a}
{\sc M.~A. M\"uller and L.~Gr\"une}, {\em Economic model predictive control
  without terminal constraints for optimal periodic behavior}, Automatica, 70
  (2016), pp.~128--139.

\bibitem{OBGF20}
{\sc R.~Ou, M.~H. Baumann, L.~Gr{\"u}ne, and T.~Faulwasser}, {\em A simulation
  study on turnpikes in stochastic {LQ} optimal control}.
\newblock arXiv Preprint 2010.12201, 2020.
\newblock To appear in the Proceedings of the 11th Symposium on Advanced
  Control of Chemical Processes, IFAC ADCHEM 2021.

\bibitem{PiAG19}
{\sc S.~Pirkelmann, D.~Angeli, and L.~Gr{\"u}ne}, {\em Approximate computation
  of storage functions for discrete-time systems using sum-of-squares
  techniques}, in Proceedings of NOLCOS 2019, vol.~52 of IFAC-PapersOnLine,
  Vienna, Austria, 2019, pp.~508--513.

\bibitem{Porr18}
{\sc A.~Porretta}, {\em On the turnpike property for mean field games}, Minimax
  Theory Appl., 3 (2018), pp.~285--312.

\bibitem{PorZ13}
{\sc A.~Porretta and E.~Zuazua}, {\em Long time versus steady state optimal
  control}, SIAM J. Control Optim., 51 (2013), pp.~4242--4273.

\bibitem{QB2003}
{\sc S.~Qin and T.~Badgwell}, {\em A survey of industrial model predictive
  control technology}, Control Engineering Practice, 11 (2003), pp.~733--764.

\bibitem{Rams28}
{\sc F.~P. Ramsey}, {\em A mathematical theory of saving}, The Economic
  Journal, 38 (1928), pp.~543--559.

\bibitem{RaMD17}
{\sc J.~B. Rawlings, D.~Q. Mayne, and M.~M. Diehl}, {\em Model Predictive
  Control: Theory, Computation and Design}, Nob Hill Publishing, Madison,
  Wisconsin, 2017.

\bibitem{RebA11}
{\sc M.~Reble and F.~Allg\"ower}, {\em Unconstrained model predictive control
  and suboptimality estimates for nonlinear continuous-time systems},
  Automatica, 48 (2011), pp.~1812--1817.

\bibitem{RBKA19}
{\sc A.~Romer, J.~Berberich, J.~K{\"o}hler, and F.~Allg{\"o}wer}, {\em One-shot
  verification of dissipativity properties from input-output data}, IEEE
  Control Systems Letters, 3 (2019), pp.~709--714.

\bibitem{SchW00}
{\sc C.~Scherer and S.~Weiland}, {\em Linear matrix inequalities in control}.
\newblock Lecture Notes, Dutch Institute for Systems and Control, Delft, The
  Netherlands, 2000.

\bibitem{SeJK97}
{\sc R.~Sepulchre, M.~Jankovic, and P.~Kokotovi\'{c}}, {\em Constructive
  Nonlinear Control}, Springer-Verlag, Berlin, 1997.

\bibitem{SSZJ92}
{\sc S.~Sethi, H.~M. Soner, Q.~Zhang, and J.~Jiang}, {\em Turnpike sets and
  their analysis in stochastic production planning problems}, Math. Oper. Res.,
  17 (1992), pp.~932--950.

\bibitem{TrZZ18}
{\sc E.~Tr\'{e}lat, C.~Zhang, and E.~Zuazua}, {\em Steady-state and periodic
  exponential turnpike property for optimal control problems in {H}ilbert
  spaces}, SIAM J. Control Optim., 56 (2018), pp.~1222--1252.

\bibitem{TreZ15}
{\sc E.~Tr{\'e}lat and E.~Zuazua}, {\em The turnpike property in
  finite-dimensional nonlinear optimal control}, J. Differ. Equ., 258 (2015),
  pp.~81--114.

\bibitem{TuMT06}
{\sc S.~E. Tuna, M.~J. Messina, and A.~R. Teel}, {\em Shorter horizons for
  model predictive control}, in Proceedings of the 2006 American Control
  Conference, Minneapolis, Minnesota, USA, 2006, pp.~863--868.

\bibitem{vNeu38}
{\sc J.~von Neumann}, {\em {\"U}ber ein \"o{}konomisches {G}leichungssystem und
  eine {V}erallgemeinerung des {B}rouwerschen {F}ixpunktsatzes}, in
  {E}rgebnisse eines {M}athematischen {S}eminars, K.~Menger, ed., 1938.

\bibitem{Will71}
{\sc J.~C. Willems}, {\em Least squares stationary optimal control and the
  algebraic {R}iccati equation}, IEEE Trans. Autom. Control, 16 (1971),
  pp.~621--634.

\bibitem{Will72a}
{\sc J.~C. Willems}, {\em Dissipative dynamical systems. {I}. {G}eneral
  theory}, Arch. Rational Mech. Anal., 45 (1972), pp.~321--351.

\bibitem{Will72b}
{\sc J.~C. Willems}, {\em Dissipative dynamical systems. {II}. {L}inear systems
  with quadratic supply rates}, Arch. Rational Mech. Anal., 45 (1972),
  pp.~352--393.

\bibitem{ZanF18}
{\sc M.~Zanon and T.~Faulwasser}, {\em Economic mpc without terminal
  constraints: Gradient-correcting end penalties enforce asymptotic stability},
  Journal of Process Control, 63 (2018), pp.~1--14.

\bibitem{ZaGD13}
{\sc M.~Zanon, S.~Gros, and M.~Diehl}, {\em A {L}yapunov function for periodic
  economic optimizing model predictive control}, in Proceedings of the 52nd
  IEEE Conference on Decision and Control --- CDC2013, Florence, Italy, 2013,
  pp.~5107--5112.

\bibitem{ZaGD14}
{\sc M.~Zanon, S.~Gros, and M.~Diehl}, {\em Indefinite linear {M}{P}{C} and
  approximated economic mpc for nonlinear systems}, Journal of Process Control,
  24 (2014), pp.~1273--1281.
\newblock Economic nonlinear model predictive control.

\bibitem{ZaGD17}
{\sc M.~Zanon, L.~Gr\"une, and M.~Diehl}, {\em Periodic optimal control,
  dissipativity and {M}{P}{C}}, IEEE Trans. Auto. Cont., 62 (2017),
  pp.~2943--2949.

\bibitem{Zuaz17}
{\sc E.~Zuazua}, {\em Large time control and turnpike properties for wave
  equations}, Ann. Rev. Control, 44 (2017), pp.~199--210.

\end{thebibliography}

%
%
\newpage

\section[Available Storage and Required Supply]{Appendix: Available Storage and Required Supply}\label{sidebar:avstorage}

Optimal control can be used in order to compute storage functions for dissipative systems, provided the supply rate $s$ and---in case of strict dissipativity---the $\KK_\infty$-function $\alpha$ are known. More precisely, the system is strictly dissipative if and only if the optimal value function 
\[  V(x_0) := \sup_{T\ge 0, u\in\LG{\UU(x_0,T)}} \int_0^T - s(x(\tau),u(\tau)) + \alpha(\|x(\tau)-x^e\|)d\tau\] 
is finite for all initial values $x_0$. In this case, $\lambda = V$ is a storage function called the \emph{available storage}. An analogous construction works without $\alpha$ in case of non-strict dissipativity and with a sum instead of the integral in case of discrete-time systems. 

That $V$ is indeed a storage function follows for any $t>0$ and $\hat u\in\UU(x_0,t)$ from the inequalities
\begin{eqnarray*}
V(x_0) & = & \sup_{T\ge 0, u\in\LG{\UU(x_0,T)}} \int_0^T - s(x(\tau),u(\tau)) + \alpha(\|x(\tau)-x^e\|)d\tau\\
& \ge & \sup_{T\ge t, u\in\LG{\UU(x_0,T)}} \int_0^T - s(x(\tau),u(\tau)) + \alpha(\|x(\tau)-x^e\|)d\tau\\
& \ge  & \int_0^t - s(x(\tau),\hat u(\tau)) + \alpha(\|x(\tau)-x^e\|)d\tau\\
&& \;\; + \;\; \sup_{T\ge 0, u\in\LG{\UU(x(t),T-t)}} \int_t^{T-t} - s(x(\tau),u(\tau)) + \alpha(\|x(\tau)-x^e\|)d\tau\\
& = & \int_0^t - s(x(\tau),\hat u(\tau)) + \alpha(\|x(\tau)-x^e\|)d\tau  + V(x(t)),
\end{eqnarray*}
which implies \eqref{eq:ctsdiss} for $\lambda=V$. 

If each $x\in \R^n$ can be reached from the equilibrium $x^e$, then another optimal control characterization of a storage function is given by the \emph{required supply}. In this case one defines
\[ V(x) := \inf_{T\ge 0, u\in\LG{\UU(x_0,T)}: \atop x(0) = x^e, x(T)=x}\int_0^T  s(x(\tau),u(\tau)) - \alpha(\|x(\tau)-x^e\|)d\tau. \] 
Then strict dissipativity holds if and only if $V$ is bounded from below and then, again, $\lambda=V$ is a storage function. Here, the storage function property follows from the fact that steering from $x^e$ to $x(t)$ via $x_0$ cannot be cheaper than steering from $x^e$ to $x(t)$ in the optimal way. For details on both constructions we refer to \cite{Will72a}. 

\newpage

\section[Model Predictive Control]{Appendix: Model Predictive Control}\label{sidebar:mpc}

Model Predictive Control (MPC) is one of the most successful optimization based control techniques with ample applications in industry \cite{QB2003,FPHG15}. We describe it here in discrete time, noting that the adaptation to continuous time is relatively straightforward. For further reading we recommend, e.g., the monographs \cite{RaMD17,GruP17}.

In many regulation problems, one would like to solve the infinite horizon optimal control problem 
\be \mbox{minimize } J_\infty(x_0,u) = \sum_{k=0}^{\infty} \ell(x(t),u(t)) \label{eq:dtoci} \ee
with respect to $u\in\UU$, $u(t)\in\U$ and $x(t)\in \X$ for all $t=0,1,2,\ldots$, where $\UU$ is an appropriate space of sequences, $\X$ and $\U$ are as above, and $x(t)$ satisfies $x(0)=0$ and \eqref{eq:dtsys}. Unless the problem has a very particular structure (as, e.g., linear dynamics, quadratic costs and no constraints), a closed-form solution of \eqref{eq:dtoci} is usually not available and due to the infinite time horizon a numerical solution is very costly to obtain, particularly if one wants to obtain the optimal control in feedback form, i.e., in the form $u^*(t) = F(x(t))$ for a suitable map $F$. 

The key idea of MPC now consists in truncating the optimization horizon and instead of \eqref{eq:dtoci} solve \eqref{eq:dtoc} or a variant thereof. This variant may include additional terminal constraints of the form $x(T)\in\X_0$ and/or a terminal cost of the form $F(x(T))$ as an additional summand in $J_T$. The MPC loop then proceeds as follows:
\begin{enumerate}
\item[1)] Pick an initial value $x_{MPC}(0)$ and a time horizon $T$. Set $k:=0$.
\item[2)] Set $x_0:=x_{MPC}(k)$ and solve \eqref{eq:dtoc} with this initial value. Denote the resulting optimal control sequence by $u^*(\cdot)$ and set $F_{MPC}(x_{MPC}(k)) := u^*(0)$.
\item[3)] Apply $F_{MPC}(x_{MPC}(k))$, measure $x_{MPC}(k+1)$, set $k:=k+1$ and go to 2). 
\end{enumerate}
The resulting solution $x^{MPC}(t)$ is called the \emph{MPC closed-loop solution} while the individual open-loop finite horizon optimal solutions computed by solving \eqref{eq:dtoc} in Step 2 are called the \emph{predictions}.

\begin{figure}[htb]
\includegraphics[width=12cm]{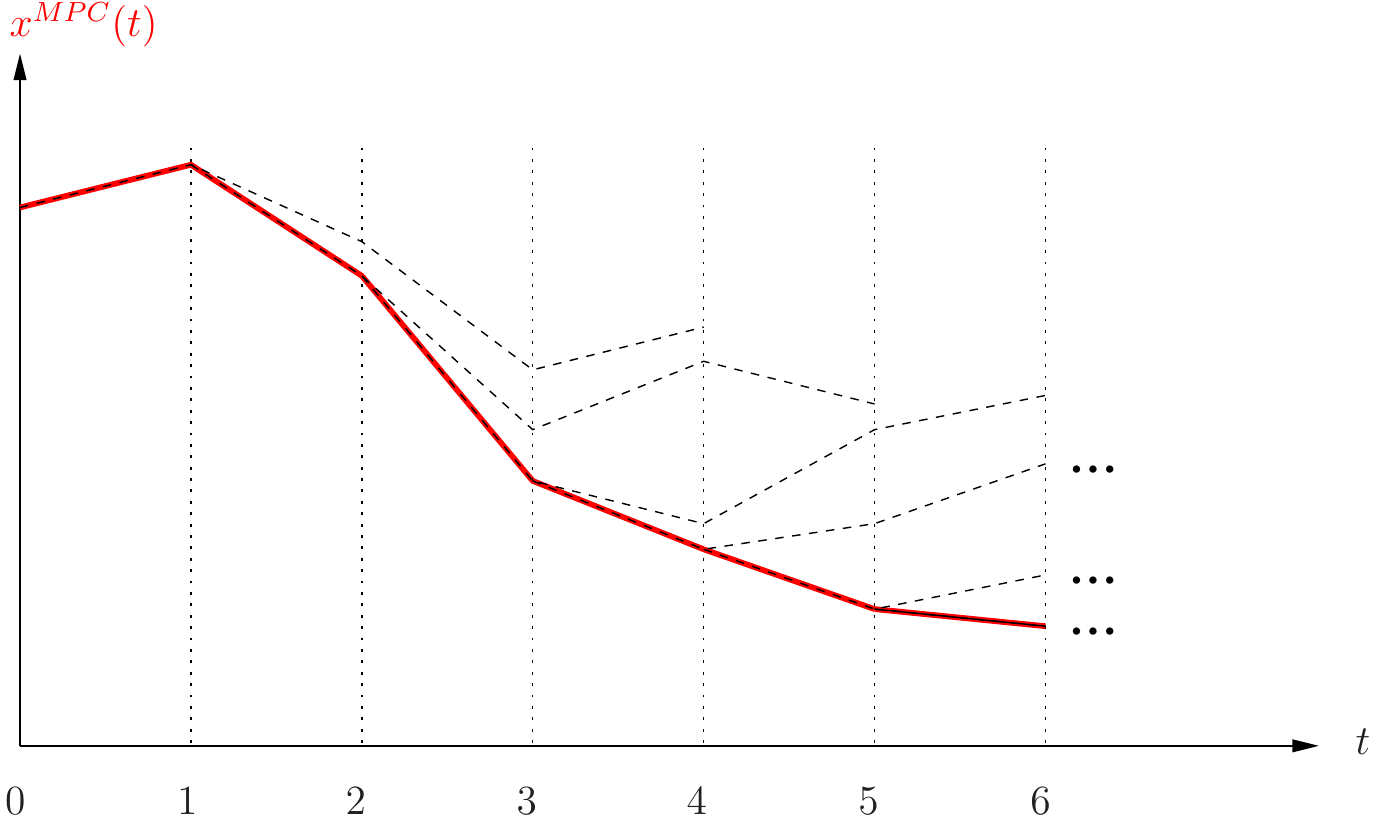}
\caption{Sketch of the solutions generated in the MPC loop. The red solid line depicts the MPC closed-loop solution while the black dashed lines indicate the predictions. The figure sketches the ideal case in which the closed-loop dynamics exactly coincides with the dynamics used to compute the predictions.\label{fig:mpcsketch}}
\end{figure}

Figure \ref{fig:mpcsketch} shows a schematic sketch of the resulting solutions loop. Here the red solid line depicts the MPC closed-loop solution while the black dashed lines indicate the predictions. The figure shows the ideal case in which the closed-loop dynamics exactly coincides with the dynamics used to compute the predictions. 

We note that by means of the dynamic programming principle (see, e.g., \cite[Theorem 4.6]{GruP17}) the solution strategy ``solve the optimal control problem and use the first element of the resulting optimal control sequence as feedback control value'' would yield an optimal feedback law if we solved the infinite horizon problem \eqref{eq:dtoci} in Step 2. However, by resorting to the (numerically much more easy to obtain) solution of the finite horizon problem \eqref{eq:dtoc} in Step 2, we make an error that needs to be analyzed. MPC for stage costs considered in this paper, which do not merely penalize the distance of the state from a desired steady state, is often denoted as \emph{economic MPC}, although the term \emph{general MPC} would probably be more fitting. In any case, strict dissipativity of the optimal control problem plays an important role in the analysis of such MPC schemes.

\end{document}